\documentclass[a4paper, 11pt, oneside]{article}

\usepackage{amsmath}
\usepackage[runin]{abstract}
\usepackage{hyperref}
\usepackage{fancyvrb}
\usepackage{amssymb}
\usepackage{mathtools}
\usepackage{amsthm}
\usepackage{amsfonts}
\usepackage{dsfont}
\usepackage{fancybox}
\usepackage{xcolor}
\usepackage{verbatim}

\usepackage{graphics,graphicx}

\newcommand{\cs}{{\mathcal S}}

\newcommand{\cf}{{\mathcal F}}
\newcommand{\cfi}{{\cf}^{-1}}

\let\svthefootnote\thefootnote
\newcommand\freefootnote[1]{%
	\let\thefootnote\relax%
	\footnotetext{#1}%
	\let\thefootnote\svthefootnote%
}

\newtheorem*{algorithm*}{Algorithm}

\theoremstyle{plain}
\newtheorem{Theorem}{Theorem}[section]
\newtheorem{Lemma}[Theorem]{Lemma}
\newtheorem{Proposition}[Theorem]{Proposition}

\newtheorem{Remark}[Theorem]{Remark}
\newtheorem{Definition}[Theorem]{Definition}
\numberwithin{equation}{section}

\newcommand{\supp}{{\rm supp \, }}
\newcommand{\spa}{{\rm span \, }}

\newtheorem{example}[Theorem]{Example}

\author{Marc Hovemann$^{1,a,\star}$, Anne Kopsch$^{2,a}$   , Thorsten Raasch$^{3,b}$, \\ Dorian Vogel$^{4,a}$ }
\title{B-Spline Quarklets and Biorthogonal Multiwavelets}
\date{\today}

\begin{document}
\maketitle
\freefootnote{This work has been supported by Deutsche Forschungsgemeinschaft (DFG), grants DA360/24-1 and RA2090/3-1.

$\star$ : Corresponding author.
	
a: Philipps-Universität Marburg, Fachbereich Mathematik und Informatik,   Hans-Meerwein Str. 6,  Lahnberge, 35043 Marburg, Germany.
	
b: Universität Siegen, Department Mathematik,  Walter-Flex-Str.
	3, 57068 Siegen, Germany.
	
E-mail: 1: \texttt{hovemann@mathematik.uni-marburg.de}, 2: \texttt{kopscha@mathematik.uni-marburg.de}, 3: \texttt{raasch@mathematik.uni-siegen.de}, 4:  
\texttt{vogeldor@mathematik.uni-marburg.de}	
	}

\textbf{Abstract.} We show that B-spline quarks and the associated quarklets fit into the theory of biorthogonal multiwavelets. Quark vectors are used to define sequences of subspaces $  V_{p,j}  $ of $ L_{2}(\mathbb{R})   $ which fulfill almost all conditions of a multiresolution analysis. Under some special conditions on the parameters they even satisfy all those properties. Moreover we prove that quarks and quarklets possess modulation matrices which fulfill the perfect reconstruction condition. Furthermore we show the existence of generalized dual quarks and quarklets which are known to be at least compactly supported tempered distributions from $\mathcal{S}'(\mathbb{R})$.  Finally we also verify that quarks and quarklets can be used to define sequences of subspaces $ W_{p,j}    $ of $  L_{2}(\mathbb{R}) $ that yield non-orthogonal decompositions of $ L_{2}(\mathbb{R})  $.

\vspace{0,3 cm}

\textbf{MSC 2020.}  42C40, 41A15.

\vspace{0,3 cm}

\textbf{Key Words.} B-spline quarks, quarklets, biorthogonal CDF-wavelets, biorthogonal multiwavelets, multiresolution analysis, perfect reconstruction condition.

\section{Introduction}

In this paper our subjects of interest are B-spline quarks and associated quarklets. We investigate how they fit into the theory of biorthogonal multiwavelets. Originally quarkonial decompositions have been introduced by Triebel in \cite{Tr01}, see also \cite{DaORaa}.  The quarks and quarklets we are concerned with in this paper have been constructed in the last decade by Dahlke, Keding and Raasch, see \cite{bib:DKR17}. Later their properties have been studied in detail in \cite{bib:DFK18}, \cite{bib:DRS19}, \cite{bib:sieber2020adaptive} and \cite{HoDa2021}. The quarks are cardinal B-splines that are multiplied with some monomial of degree $ p \in \mathbb{N}_{0}   $. The associated quarklets are constructed by means of biorthogonal compactly supported Cohen-Daubechies-Feauveau spline wavelets, where the primal generators are cardinal B-splines. For the theory of such biorthogonal wavelets we refer to Section 6.A in \cite{bib:CDF92}. Roughly speaking our quarklets $ \psi_{p}  $ are linear combinations of translated cardinal B-splines that are multiplied with some monomial of degree $ p \in \mathbb{N}_{0}  $. Precise explanations can be found in the Definitions \ref{Bquark} and \ref{def_quarklet} below. One important motivation to study quarks and quarklets are some
very interesting applications in the context of the numerical treatment of elliptic partial differential equations. Here adaptive finite element methods (AFEM) are well-established tools. The classical $h$-FEM relies on a space refinement, whereas for $p$-methods the polynomial degree of the ansatz functions varies. A combination of both methods is also possible and called $hp$-FEM. For an overview of FEM we refer to \cite{Ci}, \cite{Ha} and \cite{Sch}. There also exist adaptive numerical schemes based on wavelets that are guaranteed to converge with optimal order, see for example \cite{bib:CDD01}. Those strategies can be interpreted as $h$-methods. Therefore the question arises whether
it is possible to design $hp$-versions of adaptive wavelet schemes. This issue directly leads to our quarklets. Some first results concerning adaptive quarklet methods can be found in \cite{bib:DFK18}, \cite{bib:DKR17}, \cite{bib:DRS19} and \cite{bib:sieber2020adaptive}. In the present paper it is our main goal to discover the connections between quarklets and the theory of biorthogonal multiwavelets. Some first findings concerning this topic already have been made in \cite{bib:sieber2020adaptive} and \cite{bib:DRS19}. So it is already known that quark vectors are refinable, see Proposition \ref{prop_ref1}. Moreover it was observed that quarks and quarklets are connected via decomposition relations, see Lemma \ref{quarklet_ref1}. Using these results as a starting point in this paper we will see that our quarklets fit into the theory of biorthogonal multiwavelets very well. So (beside others) we obtained the following findings.

\begin{itemize}
\item[(I)] Quark vectors can be used to define sequences of subspaces $  V_{p,j}  $ of $ L_{2}(\mathbb{R})   $ which fulfill almost all conditions of a multiresolution analysis. Under some special conditions on the parameters they even have all properties of a multiresolution analysis. Then quark vectors can be interpreted as multiscaling functions with multiplicity $p+1$. For the details we refer to Section \ref{sec_multiscale} and Theorem \ref{thm_sec3summ1}.

\item[(II)] Quarks and quarklets possess modulation matrices that fulfill the perfect reconstruction condition. For that we refer to Section \ref{sec_perrec} and Theorem \ref{Satz_qua_perref1}. 

\item[(III)] There exist generalized dual quarks and quarklets which are at least compactly supported tempered distributions from $\mathcal{S}'(\mathbb{R})$. Here generalized means that due to the perfect reconstruction property there exists a dual modulation matrix which implies a vector refinement equation that has a distributional solution. We also can give a formula to describe the generalized duals. Moreover under some very special conditions on the parameters it is possible to prove additional smoothness and integrability properties for the duals. More information can be found in Section \ref{sec_dualqua} and Lemma \ref{lem_qua_dualconst1}. 

\item[(IV)] Quarks and quarklets can be used to define sequences of subspaces $ W_{p,j}    $ of $  L_{2}(\mathbb{R}) $ that yield non-orthogonal decompositions of $ L_{2}(\mathbb{R})  $. It is also possible to rewrite the spaces $ W_{p,j}    $ by using modified quarklets with additional orthogonality properties to obtain an alternative decomposition of $  L_{2}(\mathbb{R})   $. For that we refer to Section \ref{sec_decompoqua}, see Theorem \ref{satz_L2zerl1} and Theorem \ref{satz_L2zerl2}.  

\end{itemize}

The authors are confident that the findings from this paper also can be used to push forward the research concerning adaptive quarklet approximation. Some of the results of this article seem to allow the conclusion that when it comes to practical applications it is better to work with quarklets that are based on B-splines with odd $m \in \mathbb{N}$. Here $m$ is as in \eqref{eq_defm}. Some findings that substantiate this hypothesis can be found in Lemma \ref{lem_stab_sinleres}. Another outcome of this paper which seems to be valuable for the numerical application of the quarklets is the discovery of the modified quarklets $ \Psi^{\star}  $, see Definition \ref{def_Wpj_star}. Due to their additional orthogonality properties they already have been used to improve the results concerning adaptive quarklet tree approximation at least in the case of Haar quarklets, namely $ m=1$, see \cite{bib:DaHoRaaVo}.  

This paper is organized as follows. In Section \ref{sec_quarkdef} we give a precise definition of the quarks and quarklets. In Section \ref{sec_multiscale} we show how the quarks can be used to construct a multiresolution analysis. Section \ref{sec_perrec} explains how quarks and quarklets fit into the theory of biorthogonal multiwavelets. Here we prove that quarks and quarklets fulfill the perfect reconstruction condition. In the course of Section \ref{sec_dualqua} we prove the existence of generalized dual quarks and quarklets. Moreover we collect some elementary properties of them. Finally in Section \ref{sec_decompoqua} we illustrate how quarklets can be used to find a (non-orthogonal) decomposition of $  L_{2}(\mathbb{R})$. Furthermore we describe a method to orthogonalize the quarklets. 

First of all we fix some notation. As usual $\mathbb{N}$ denotes the natural numbers, $\mathbb{N}_0$ the natural numbers including $0$, $\mathbb{Z}$ the integers, $\mathbb{R}$ the real numbers and $ \mathbb{C}  $ the complex numbers.  Let $\mathcal{S}(\mathbb{R})$ be the collection of all Schwartz functions on $\mathbb{R}$ endowed with the usual topology and denote by $\mathcal{S}'(\mathbb{R})$ its topological dual, namely the space of all bounded linear functionals on $\mathcal{S}(\mathbb{R})$ endowed with the weak $\ast$-topology. For $ 0 < r \leq \infty $ by $ L_{r}(\mathbb{R})  $ we denote the usual Lebesgue spaces. Given a function $ f \in  L_{r}(\mathbb{R}) $ we use the symbol $ \Vert f \vert L_{r}(\mathbb{R}) \Vert   $ for the associated quasi-norm. For $ f \in L_{1}(\mathbb{R})  $ we define its Fourier transform via
\begin{align*}
( \mathcal{F}f)(\xi) := \frac{1}{\sqrt{2 \pi}} \int_{\mathbb{R}} f(x) e^{- i x \xi} dx .
\end{align*} 
The symbol $ \mathcal{F}^{-1}  $ refers to the inverse transform. Recall that both $\cf$ and $\cfi$ also can be defined on $\cs'(\mathbb{R})$. For $ f, g \in L_{2}(\mathbb{R})  $ we define the inner product
\begin{align*}
\left\langle f , g     \right\rangle_{L_{2}(\mathbb{R})} := \int_{\mathbb{R}} f(x) \overline{g(x)} dx .
\end{align*}
Given a set of functions $G \subset L_{2}(\mathbb{R})$ by $ \overline{G} $ we mean the $ L_{2}-$closure of $G$. For a finite set of functions $  f_{1}, \ldots , f_{n} \in L_{2}(\mathbb{R})    $ their closed linear span 
\begin{align*}
V = \overline{ \spa \{ f_{l}( \cdot - k ) : 1 \leq l \leq n , k \in \mathbb{Z}   \} } \subset L_{2}(\mathbb{R}) 
\end{align*}
is defined by
\begin{equation}\label{eq_def_span_colsedL2}
V = \overline{  \Big \{ \sum_{l = 1}^{n} \sum_{k \in \mathbb{Z}} a_{l}(k) f_{l}( \cdot - k ) :  a_{l} \in l_{0}(\mathbb{Z})  \Big \} } .
\end{equation}
Thereby $  l_{0}(\mathbb{Z})  $ is the space of all sequences that have only finitely many nonzero elements. If a function $f$ is at least $n$-times differentiable, then we denote the $n$-th derivative by $ f^{(n)}  $.  The symbols  $C, C_1, c, c_{1} \ldots $ denote constants that depend only on the fixed parameters. Unless otherwise stated their values may vary from line to line. When we write $a\sim b$ we mean that there are constants $0 < C_1 \le C_2 < \infty$ such that $a \le C_{1} b \le C_2 a$. 

\section{Quarks and Quarklets}\label{sec_quarkdef}
In this section we present a precise definition of quarks and quarklets. For that purpose in a first step we recall the definition of cardinal B-splines. Cardinal B-splines are defined recursively by $ N_{1} := \chi_{[0,1)} $ and for $  m \in \mathbb{N}   $ with $ m > 1  $ by
\begin{equation}\label{eq_defm}
N_{m} := N_{m - 1} \ast N_{1} = \int_{0}^{1} N_{m - 1}(\cdot - t) dt .
\end{equation}
Those functions possess some very nice properties. They are collected in the following lemma.

\begin{Lemma}\label{Bspl_elem}
Let $ m \in \mathbb{N}   $ and $ x \in \mathbb{R}  $. Then for the cardinal B-splines the following properties hold. 
	
	\begin{itemize}
		\item[(i)] We can write $  N_{m}(x) = \frac{1}{(m-1)!} \sum_{k = 0}^{m}(-1)^{k} { m \choose k } \max\{0,(x - k)\}^{m-1}   $.
		
		\item[(ii)] For $ m \geq 2  $ the recursion formula $ N_{m}(x) = \frac{x}{m-1} N_{m-1}(x) + \frac{m-x}{m-1} N_{m-1} (x-1)   $ holds.
		
		\item[(iii)] For $ m \geq 3  $ the derivatives satisfy $  N_{m}'(x) = N_{m-1}(x) - N_{m-1}(x-1)  $.
		
		\item[(iv)] The B-splines are compactly supported with $ \supp N_{m} = [0,m]  $.
		
		\item[(v)] We have $ \sum_{k \in \mathbb{Z} } N_{m}(x - k) = 1   $.
		
		\item[(vi)] The cardinal B-splines are refinable with $ N_{m}(x) = \sum_{k = 0}^{m} 2^{1-m} { m \choose k } N_{m}(2x - k) $.
		
	\end{itemize}
\end{Lemma}
Those properties are well-known. Let us refer to \cite{bib:devore1993constructive}, see Chapters 5.2 and 5.3.  One can also consult \cite{bib:Chui92} and Section IX in \cite{CBoor1}. In what follows for fixed $ m \in \mathbb{N} $ we will work with the symmetrized cardinal B-spline $ \varphi(x)  :=   N_{m}  ( x + \lfloor \frac{m}{2} \rfloor   )    $. We observe $ \supp \varphi = [ - \lfloor \frac{m}{2} \rfloor    ,  \lceil \frac{m}{2} \rceil    ]   $. The symmetrized cardinal B-spline shows up in the following definition where we explain the so-called quarks.

\begin{Definition}\label{Bquark}
	Let $ m \in \mathbb{N} $ and $ p \in \mathbb{N}_{0}  $. Then the \emph{$p$-th cardinal B-spline quark} $ \varphi_{p}  $ is defined by
	\begin{equation}
	\varphi_{p}(x)  := \Big ( \frac{x}{\lceil \frac{m}{2} \rceil } \Big )^{p}  N_{m} \Big ( x + \lfloor \frac{m}{2} \rfloor  \Big ) .
	\end{equation}
\end{Definition}
The quarks will be very important for us in order to define the quarklets. Their properties have been studied in \cite{bib:DKR17}.

It is shown in \cite{bib:CDF92} by Cohen, Daubechies and Feauveau, that for a given $ \tilde{m} \in \mathbb{N}  $ with $ \tilde{m} \geq m    $ and $  m + \tilde{m} \in 2 \mathbb{N}   $ there exists a compactly supported spline wavelet $ \psi  $ (sometimes called CDF-wavelet) with
\begin{equation}\label{def_CDF_wav}
\psi = \sum_{k \in \mathbb{Z}} b_{k}  \varphi ( 2 \cdot - k   )
\end{equation}
with expansion coefficients $ b_{k} \in \mathbb{R}   $. Only finitely many of them are not zero. Moreover $ \psi  $ has $ \tilde{m}   $ vanishing moments and the  system
\begin{align*}
\Big \{ \varphi (  \cdot - k )  \ : \ k \in \mathbb{Z}  \Big  \} \cup \Big \{ 2^{\frac{j}{2}} \psi (2^{j} \cdot - k) \ : \ j \in \mathbb{N}_{0} \ , \ k \in \mathbb{Z}  \Big \}
\end{align*}
is a Riesz basis for $ L_{2}(\mathbb{R})   $. To construct such a $ \psi   $ we have to work with a compactly supported dual generator $ \tilde{\varphi}   $ associated to the primal generator $ \varphi $ that fulfills
\begin{equation}\label{biorto1}
\left \langle  \varphi , \tilde{\varphi} (\cdot - k)   \right\rangle_{L_{2}(\mathbb{R})} = \delta_{0,k} , \qquad k \in \mathbb{Z}.
\end{equation}
Connected with that it is possible to construct another compactly supported wavelet $ \tilde{\psi} \in L_{2}(\mathbb{R})  $ with 
\begin{equation}\label{def_biort_wav1}
\tilde{\psi} = \sum_{k \in \mathbb{Z}}  \tilde{b}_{k} \tilde{\varphi} ( 2 \cdot - k   ).
\end{equation}
Here only finitely many of the $ \tilde{b}_{k} \in \mathbb{R}   $  are not zero. Moreover $ \tilde{\psi}  $ has $ m \in \mathbb{N}  $ vanishing moments. For $ j \in \mathbb{N}_{0}   $ and $ k \in \mathbb{Z}   $ let us write
\begin{equation}\label{def_biort_wav2}
\psi_{j,k} = 2^{\frac{j}{2}} \psi ( 2^{j} \cdot - k ) \qquad \qquad \mbox{and} \qquad \qquad \tilde{\psi}_{j,k} = 2^{\frac{j}{2}} \tilde{\psi} ( 2^{j} \cdot - k ) .
\end{equation}
For $ k \in \mathbb{Z}   $ we put $ \psi_{-1,k} = \varphi (  \cdot - k )   $ and $ \tilde{\psi}_{-1,k} = \tilde{\varphi} (  \cdot - k )   $. Then both function systems are biorthogonal in the following sense 
\begin{equation}\label{biorto2}
\langle \psi_{j,k} , \tilde{\psi}_{j',k'}     \rangle_{L_{2}(\mathbb{R})} = \delta_{j , j'} \delta_{k , k'} , \qquad j, j' \in \mathbb{N}_{0}, \quad k, k' \in \mathbb{Z} .
\end{equation}
Furthermore for each $ f \in L_{2}(\mathbb{R})  $ we have
\begin{align}
f & = \sum_{k \in \mathbb{Z}} \langle f , \tilde{\psi}_{-1,k}     \rangle_{L_{2}(\mathbb{R})} \psi_{-1,k}  + \sum_{j \in \mathbb{N}_{0} ,k \in \mathbb{Z}} \langle f , \tilde{\psi}_{j,k}     \rangle_{L_{2}(\mathbb{R})} \psi_{j,k} \nonumber \\
& = \sum_{k \in \mathbb{Z}} \langle f , \psi_{-1,k}     \rangle_{L_{2}(\mathbb{R})} \tilde{\psi}_{-1,k}  +  \sum_{j \in \mathbb{N}_{0}, k \in \mathbb{Z}} \langle f , \psi_{j,k}     \rangle_{L_{2}(\mathbb{R})} \tilde{\psi}_{j,k}    \label{def_biort_wav3}
\end{align}
with convergence in $ L_{2}(\mathbb{R})  $. For details and proofs concerning the above construction we refer to Section 6.A in \cite{bib:CDF92}. Now we can use the CDF-wavelets $ \psi $ to define the quarklets.  
\begin{Definition}\label{def_quarklet}
Let $ m \in \mathbb{N}  $ and $ \tilde{m} \in \mathbb{N}    $ with $ m \leq \tilde{m}   $ and $ m + \tilde{m} \in 2 \mathbb{N}   $. 	Let $ p \in \mathbb{N}_{0}  $. Then the \emph{$p$-th quarklet} $ \psi_{p} $ is defined by 
	\begin{equation}
	\psi_{p} := \sum_{k \in \mathbb{Z}} b_{k} \varphi_{p}(2 \cdot - k).
	\end{equation}
	Here the $ b_{k}  $ are the same as in \eqref{def_CDF_wav}. 
\end{Definition}
It can be shown that the quarks and quarklets inherit some important properties of the B-splines and B-spline wavelets, respectively. In particular Jackson and Bernstein estimates can be proved and the quarklets possess the same amount of vanishing moments. For details concerning this topic we refer to  \cite{bib:DKR17}.

\section{Quarks as Multigenerators:  Multiscale Approximation}\label{sec_multiscale}

It is one of the main goals of this paper to investigate how the quarks and quarklets are connected with the theory of biorthogonal multiwavelets. For that purpose in a first step we deal with the quark vector 
\begin{equation}\label{eq_quark_vec1}
\Phi :=  ( \varphi_{0} , \varphi_{1} , \ldots , \varphi_{p}    )^{T}        
\end{equation}
with $ p \in \mathbb{N}_{0}   $, where the quarks are as in Definition \ref{Bquark}. We want to find out whether $ \Phi $ can be used to construct a multiresolution analysis. Let us start with recalling the definition of $L_{2}-$stability for function vectors. 

\begin{Definition}\label{def_L2stab}
Let $ F = ( f_{1} , \ldots , f_{n}    )^{T} \in (L_{2}(\mathbb{R}))^{n}   $ be a function vector. Then we say that $   F  $ is \emph{$L_{2}-$stable}, if there are constants $  0 < C_{1} \leq C_{2} < \infty   $ such that
\begin{align*}
C_{1} \sum_{k = - \infty}^{\infty} \overline{v}^{T}_{k} v_{k} \leq \Big \Vert \sum_{k = - \infty}^{\infty} \overline{v}^{T}_{k} F (\cdot - k) \Big \vert L_{2}(\mathbb{R}) \Big \Vert^{2} \leq C_{2} \sum_{k = - \infty}^{\infty} \overline{v}^{T}_{k} v_{k}
\end{align*}
holds for any vector sequence $  \{ v_{k}  \}_{k \in \mathbb{Z}} \in l_{2}(\mathbb{Z})^{n}   $. Here $ l_{2}(\mathbb{Z})^{n}   $ denotes the set of sequences of vectors $ v_{k} \in \mathbb{C}^{n}    $ with $  \sum_{k = - \infty}^{\infty} \overline{v}^{T}_{k} v_{k} < \infty   $.
\end{Definition}

There exist several possibilities to equivalently describe the $L_{2}-$stability from Definition \ref{def_L2stab}. One of them is given in the following result. 

\begin{Lemma}\label{lem_L2stab}
Let $ f_{1} , \ldots , f_{n}  \in L_{2}(\mathbb{R})     $ be compactly supported. Then the integer translates of $ f_{1} , \ldots , f_{n}      $ are $L_{2}-$stable if and only if for any $ \xi \in \mathbb{R}   $ the sequences $ ( (\mathcal{F} f_{l} )(\xi + 2 \pi \alpha)     )_{\alpha \in \mathbb{Z}}     $ with $  l = 1 , \ldots , n   $ are linearly independent. 
\end{Lemma}

This result can be found in \cite{JiMi}, see Theorem 4.1. Let us also refer to Section 1.2 in \cite{PlStr}. When we deal with a single function Lemma \ref{lem_L2stab} becomes easier. Indeed in that case we observe the following.

\begin{Lemma}\label{lem_stab_singl}
Let $ f  \in L_{2}(\mathbb{R})     $ be compactly supported. Then the integer translates of $ f $ are $L_{2}-$stable if and only if for all $  \xi \in \mathbb{R}  $  we have
\begin{align*}
\sum_{\alpha \in \mathbb{Z}} | (\mathcal{F} f ) (\xi + 2 \pi \alpha )    |^{2} > 0 .
\end{align*}  
\end{Lemma}

This result can be found in \cite{JiMi}, see Theorem 3.3. Let us mention that if $ f  \in L_{2}(\mathbb{R})     $ is compactly supported, then we have $ f  \in \mathcal{L}_{2}(\mathbb{R})     $, see the remark before Theorem 2.1 in \cite{JiMi} for a definition and explanations. Now we are well prepared to recall what a multiresolution analysis in the context of multiwavelets is, see also Definition 6.10 in \cite{Kein1}. 

\begin{Definition}\label{def_MRA}
Let $ \{ V_{j}  \}_{j \in \mathbb{Z}}  $ be a sequence of closed subspaces of $  L_{2}(\mathbb{R}) $. Then this sequence is called a \emph{multiresolution analysis (MRA)} if the following properties are fulfilled.

\begin{itemize}
\item[(i)] For all $ j \in \mathbb{Z}  $ we have $ V_{j} \subset V_{j+1}  $.
\item[(ii)] It is $ \overline{ \bigcup_{j \in \mathbb{Z}} V_{j} } = L_{2}(\mathbb{R})  $.
\item[(iii)] It is $  \bigcap_{j \in \mathbb{Z}} V_{j}  = \{ 0 \}  $.
\item[(iv)] We have $ f \in V_{j}  $ if and only if $ f(\cdot - 2^{-j} k) \in V_{j}  $ for all $ j, k \in \mathbb{Z}  $.
\item[(v)] We have $ f \in V_{j}  $ if and only if $ f( 2 \cdot ) \in V_{j+1}  $ for all $ j \in \mathbb{Z}  $.
\item[(vi)] For some $ r \in \mathbb{N}  $ there exists an $  L_{2}-$stable function vector $ G = ( g_{1} , \ldots , g_{r} )^{T}   $ that consists of $r$ functions of $ L_{2}(\mathbb{R})  $ such that $ V_{0} = \overline{ \spa \{ g_{l} (\cdot - k), l \in \{ 1 , \ldots , r \} , k \in \mathbb{Z}   \} }    $. The vector $G$ is called \emph{multiscaling function} of multiplicity $r$.  
\end{itemize}
\end{Definition}

To incorporate the quarks into the setting of multiresolution analyses let us define the following subspaces of $ L_{2}(\mathbb{R})  $. Here we follow \cite{bib:DKR17} and \cite{bib:sieber2020adaptive}.

\begin{Definition}\label{def_Vpj}
Let $ j \in \mathbb{Z}  $ and $  p \in \mathbb{N}_{0}  $. Then we define the spaces
\begin{align*}
V_{p,j} :=  \overline{ \spa \{ 2^{\frac{j}{2}} \varphi_{q} (2^{j} \cdot - k ) : 0 \leq q \leq p , k \in \mathbb{Z}   \} } \subset L_{2}(\mathbb{R}).
\end{align*}
Here the functions $ \varphi_{q}  $ are the quarks from Definition \ref{Bquark}.
\end{Definition}

\begin{Remark}\label{rem_V0_multana}
It is well-known that the sequence $ \{ V_{0,j}  \}_{j \in \mathbb{Z}}  $ is a multiresolution analysis in the sense of Definition \ref{def_MRA} with multiplicity $r=1$. For that we refer to \cite{bib:CDF92}, see in particular Section 6.A.
\end{Remark}
Now let $  p  \in \mathbb{N}_{0} $. In what follows we want to investigate whether the sequences $  \{ V_{p,j}  \}_{j \in \mathbb{Z}}    $ form multiresolution analyses in the sense of Definition \ref{def_MRA}. For that purpose we examine whether the properties listed in Definition \ref{def_MRA} are fulfilled step by step. Initially let us recall that the cardinal B-spline quarks are refinable in the following sense. 

\begin{Proposition}\label{prop_ref1}
Let $ m \in \mathbb{N}  $ and $ p \in \mathbb{N}_{0}   $. Then the quark vector $  ( \varphi_{0} , \varphi_{1} , \ldots , \varphi_{p}  )^{T}  $ is refinable. More precisely there exist $ ( p + 1 ) \times ( p + 1 )$ refinement matrices $ A_{k} $ such that for all $  x \in \mathbb{R} $ the two-scale matrix refinement equation
\begin{equation}\label{eq_ref_1aaaa}
\left ( \begin{array}{c} \varphi_{0}(x) \\ \vdots  \\ \varphi_{p}(x)  \end{array}   \right ) = \sum_{k \in \mathbb{Z}} A_{k} \left ( \begin{array}{c} \varphi_{0}(2x - k) \\ \vdots  \\ \varphi_{p}(2x - k)  \end{array}   \right ) 
\end{equation} 
holds. The matrices $ A_{k} $ are given by 
\begin{align*}
( A_{k} )_{q,l} := \frac{1}{2^{q-1}} \Big ( \lceil \frac{m}{2} \rceil \Big )^{l-q}   a_{k} { q - 1 \choose  l - 1 } k^{q - l}   \qquad \qquad \mbox{with} \ q, l \in \{ 1 , \ldots , p+1 \} .
\end{align*}
Here for $  k \in \mathbb{Z}  $ with  $ - \lfloor \frac{m}{2} \rfloor \leq k \leq \lceil \frac{m}{2} \rceil  $ it is 
\begin{align*}
a_{k} = 2^{1 - m} { m \choose k +  \lfloor \frac{m}{2} \rfloor    } .
\end{align*}
\end{Proposition}

This result can be found in \cite{bib:DKR17}, see Proposition 5. Now we are well-prepared to formulate the first main result of this section. It tells us that the spaces $ V_{p,j}  $ from Definition \ref{def_Vpj} fulfill most of the properties which can be found in Definition \ref{def_MRA}.

\begin{Theorem}\label{thm_sec3summ}
Let $ m \in \mathbb{N} $ and $ p \in \mathbb{N}_{0} $. Then the sequence $  \{ V_{p,j}  \}_{j \in \mathbb{Z}}    $ of spaces from Definition \ref{def_Vpj} fulfills the properties (i)-(v) from Definition \ref{def_MRA}.
\end{Theorem}

\begin{proof}
\textit{Step 1.} At first we prove that property (i) from Definition \ref{def_MRA} is fulfilled. More precisely we show that for all $ j \in \mathbb{Z}  $ we have $   V_{p,j} \subset V_{p,j+1}  $. However this is a simple consequence of Proposition \ref{prop_ref1} and follows from formula \eqref{eq_ref_1aaaa} by standard arguments. 

\textit{Step 2.} Next we verify that the properties (iv) and (v) from Definition \ref{def_MRA} are fulfilled. Obviously they are direct consequences of the definition of the spaces $ V_{p,j} $, see Definition \ref{def_Vpj}.

\textit{Step 3.} Now we prove that property (ii) from Definition \ref{def_MRA} is fulfilled. With other words we show that it holds $ \overline{ \bigcup_{j \in \mathbb{Z}} V_{p,j} } = L_{2}(\mathbb{R})     $. The proof is based on the fact that $ \{ V_{0,j}  \}_{j \in \mathbb{Z}}  $ is a multiresolution analysis, see Remark \ref{rem_V0_multana}.  Let us start with the obvious observation $  \overline{ \bigcup_{j \in \mathbb{Z}} V_{p,j} } \subset L_{2}(\mathbb{R})      $.  Therefore it remains to prove $  \overline{ \bigcup_{j \in \mathbb{Z}} V_{p,j} } \supset L_{2}(\mathbb{R})      $.  Definition \ref{def_Vpj} implies $  \overline{ \bigcup_{j \in \mathbb{Z}} V_{p,j} } \supset  \overline{ \bigcup_{j \in \mathbb{Z}} V_{0,j} }    $. Moreover we know that the sequence $ \{ V_{0,j}  \}_{j \in \mathbb{Z}}  $ is a multiresolution analysis, see Remark \ref{rem_V0_multana}. Consequently we get
\begin{align*}
\overline{ \bigcup_{j \in \mathbb{Z}} V_{p,j} } \supset  \overline{ \bigcup_{j \in \mathbb{Z}} V_{0,j} } = L_{2}(\mathbb{R}) .
\end{align*}  
Hence this step of the proof is complete. Let us remark that property (ii) from Definition \ref{def_MRA} also can be proved by using Theorem 11.22 from \cite{Kein1}.

\textit{Step 4.} To complete the proof we show that property (iii) from Definition \ref{def_MRA} holds. That means we verify $  \bigcap_{j \in \mathbb{Z}} V_{p,j}  =  \{ 0 \}     $. For that purpose we observe that our spaces $ V_{p,j}  $ fit into the theory developed in \cite{RoZu}. So the desired result follows directly from Theorem 2.2 in \cite{RoZu}. The spaces $  V_{p,0}  $ are the smallest closed shift-invariant subspaces of $  L_{2}(\mathbb{R})   $ that are containing $ \Phi =  ( \varphi_{0} , \varphi_{1} , \ldots , \varphi_{p}    )^{T}    $, see Definition \ref{def_Vpj}. Moreover $ \Phi =  ( \varphi_{0} , \varphi_{1} , \ldots , \varphi_{p}    )^{T}    $ is a finite set of functions from $  L_{2}(\mathbb{R})   $ that generates $  V_{p,0}  $. Furthermore the spaces $ V_{p,j}  $ are the $ 2^{j}-$dilates of $ V_{p,0}  $. Consequently Theorem 2.2 from \cite{RoZu} can be applied and the result follows. Finally let us mention that there exists an alternative method to prove property (iii) from Definition \ref{def_MRA} which is based on Theorem 11.23 from \cite{Kein1}. 
\end{proof}

It remains to investigate whether the spaces $ V_{p,j}  $ satisfy property (vi) from Definition  \ref{def_MRA}. For that purpose we have to examine the $ L_{2}-$stability of the quark vectors $ \Phi  $. To this end taking into account Lemma \ref{lem_L2stab} and Lemma \ref{lem_stab_singl} we require the Fourier transforms of the quarks.    
 
\begin{Lemma}\label{lem_qua_fourtrans}
Let $ m \in \mathbb{N} $ and $ l \in \mathbb{N}_{0} $. Then for all $ \xi \in \mathbb{R}    $ the Fourier transforms of the quarks are given by
\begin{align*}
( \mathcal{F} \varphi_{l} ) (\xi)   & =  i^{l} (2 \pi)^{\frac{m}{2}} \Big ( e^{i  ( \lfloor \frac{m}{2} \rfloor - \frac{m}{2} ) \xi }  \Big ( \frac{\sin(\frac{\xi}{2})}{\frac{\xi}{2}} \Big ) ^{m}   \Big )^{(l)}  .
\end{align*}
\end{Lemma}

\begin{proof}
To prove this result we calculate
\begin{align*}
( \mathcal{F} \varphi_{l} ) (\xi)   & = \mathcal{F} \Big (  \Big ( \frac{x}{\lceil \frac{m}{2} \rceil } \Big )^{l}  N_{m} \Big ( x + \lfloor \frac{m}{2} \rfloor  \Big ) \Big ) (\xi) \\
& = \frac{1}{  \lceil \frac{m}{2} \rceil ^{l} } \mathcal{F} \Big (  x^{l}  N_{m} \Big ( x + \lfloor \frac{m}{2} \rfloor  \Big ) \Big ) (\xi) \\
& = i^{l} \frac{1}{  \lceil \frac{m}{2} \rceil ^{l} }  \Big (  \mathcal{F} \Big [ N_{m} \Big ( x + \lfloor \frac{m}{2} \rfloor  \Big ) \Big ] (\xi) \Big )^{(l)}  \\
& = i^{l} \frac{1}{  \lceil \frac{m}{2} \rceil ^{l} }  \Big ( e^{i  \lfloor \frac{m}{2} \rfloor \xi } ( \mathcal{F}  N_{m} )(\xi)  \Big )^{(l)}  \\
& = i^{l} (2 \pi)^{\frac{m}{2}}  \frac{1}{  \lceil \frac{m}{2} \rceil ^{l} }  \Big ( e^{i  \lfloor \frac{m}{2} \rfloor \xi } [ ( \mathcal{F}  N_{1} )(\xi) ]^{m}   \Big )^{(l)}  .
\end{align*}
The Fourier transform of $   N_{1} := \chi_{[0,1)}  $ is well-known. We find
\begin{align*}
( \mathcal{F} \varphi_{l} ) (\xi)   & = i^{l} (2 \pi)^{\frac{m}{2}}  \Big ( e^{i  ( \lfloor \frac{m}{2} \rfloor - \frac{m}{2} ) \xi }  \Big ( \frac{\sin(\frac{\xi}{2})}{\frac{\xi}{2}} \Big ) ^{m}   \Big )^{(l)}  .
\end{align*}
The proof is complete. 
\end{proof}

First of all we study the $L_{2}-$stability of individual quarks. Then we observe the following.

\begin{Lemma}\label{lem_stab_sinleres}
Let $  m \in \mathbb{N}  $ and $   p \in \mathbb{N}_{0} $. The quarks $ \varphi_{p}  $ are as in Definition  \ref{Bquark}.

\begin{itemize}
\item[(i)]
Then concerning $L_{2}-$stability of the integer translates we observe the following.

\vspace{0,4 cm}

\begin{tabular}[h]{c | c | c | c | c }
  & p = 0 & p = 1 & p = 2 & p = 3 \\
  \hline
  m = 1 & stable    & stable    & stable    & stable   \\
  \hline
  m = 2 & stable    & unstable    & stable    & stable   \\
  \hline
  m = 3 & stable    & stable    & stable    & stable   \\
  \hline
  m = 4 & stable    & unstable     & stable    & unstable   \\
\end{tabular}

\vspace{0,4 cm}

\item[(ii)] Let $ m \in  \mathbb{N}   $  and $ p = 0 $. Then the integer translates of $  \varphi_{0} $ are $L_{2}-$stable.

\item[(iii)] Let $ m \in  \mathbb{N}   $ be odd and $ p = 1 $. Then the integer translates of $  \varphi_{1} $ are $L_{2}-$stable.

\item[(iv)] Let $ m \in 2 \mathbb{N}  $ and $ p = 1 $. Then the integer translates of $  \varphi_{1} $ are not $L_{2}-$stable.

\end{itemize}

\end{Lemma}

This result can be proved by using Lemma \ref{lem_qua_fourtrans} and Lemma \ref{lem_stab_singl}. With other words we have to compute the zeros of the Fourier transforms of the quarks. Since this business is quite technical we postpone the details and present them in the Appendix Section \ref{sec_appendix}. Now we continue by investigating the $L_{2}-$stability of quark vectors $ \Phi =  ( \varphi_{0} , \varphi_{1} , \ldots , \varphi_{p}    )^{T}    $.

\begin{Lemma}\label{lem_L2stab_multiquark}
Let $ p \in \mathbb{N}_{0}   $ and $ m \in \mathbb{N}  $. Let $ \Phi =  ( \varphi_{0} , \varphi_{1} , \ldots , \varphi_{p}    )^{T}    $.
\begin{itemize}
\item[(i)] Let $ p=0$. Then $ \Phi $ is $L_{2}-$stable.
\item[(ii)] Let $ m=1 $ (Haar case) and $ p \in \mathbb{N}_{0}     $. Then $ \Phi $ is $L_{2}-$stable.
\item[(iii)] Let $ m > 1   $ and $ p>0$. Then $ \Phi $ is not $L_{2}-$stable.
\end{itemize}
\end{Lemma}

\begin{proof}
To prove this result we use Lemma \ref{lem_L2stab}. Moreover we apply some findings from Lemma \ref{lem_stab_sinleres}.

\textit{Step 1.} At first we prove (i). Since here we only work with a single function this directly follows from Lemma \ref{lem_stab_sinleres}.

\textit{Step 2.} Now we prove (ii). For that purpose we use Theorem 5.1 from \cite{JiMi}. Let $ l(\mathbb{Z}) $ be the linear space of all sequences. For all sequences $  \lambda^{0} , \ldots , \lambda^{p} \in  l(\mathbb{Z})    $ we define the linear mapping
\begin{align*}
L_{\varphi_{0} , \varphi_{1} , \ldots , \varphi_{p}} : ( \lambda^{0} , \ldots , \lambda^{p} ) \mapsto \sum_{l = 0}^{p} \sum_{k \in \mathbb{Z}} \lambda^{l}_{k} \varphi_{l} ( \cdot - k ) .
\end{align*} 
Theorem 5.1 in \cite{JiMi} implies that if $ L_{\varphi_{0} , \varphi_{1} , \ldots , \varphi_{p}}    $ is injective, then for all $ \xi \in \mathbb{R}   $ the sequences $ ( (\mathcal{F} \varphi_{l} )(\xi + 2 \pi \alpha)     )_{\alpha \in \mathbb{Z}}     $ with $  l = 0 , \ldots , p   $ are linearly independent. Then the $L_{2}-$stability follows from Lemma \ref{lem_L2stab}. Now because of $ m=1 $ we have $ \varphi_{l} (x) = x^{l} \chi_{[0,1)}(x)     $. Therefore we also find $ \varphi_{l} ( x - k ) = ( x - k )^{l} \chi_{[k,k+1)}(x)   $. Consequently when we want to know whether $ L_{\varphi_{0} , \varphi_{1} , \ldots , \varphi_{p}}    $ is injective, it is enough to look at a fixed $k \in \mathbb{Z}$. For example we can work with $ k=0$. Let $ p \in \mathbb{N}$. The mapping 
\begin{align*}
L^{*}_{\varphi_{0} ,  \ldots , \varphi_{p}} : ( \lambda^{0}_{0} , \ldots , \lambda^{p}_{0} ) \mapsto \sum_{l = 0}^{p}  \lambda^{l}_{0} x^{l} \chi_{[0,1)}(x)   
\end{align*}
is injective because of the different polynomial degrees for each $  l \in \{0 , \ldots , p \}   $. Therefore we can use Theorem 5.1 in \cite{JiMi} and Lemma \ref{lem_L2stab} to complete this step of the proof.

\textit{Step 3.} Now we prove (iii) for even $m$. Here we use Lemma \ref{lem_L2stab}. We show that there exists  $ \xi \in \mathbb{R}   $ such that the sequences $ ( (\mathcal{F} \varphi_{l} )(\xi + 2 \pi \alpha)     )_{\alpha \in \mathbb{Z}}     $ with $  l = 0 , \ldots , p   $ are not linearly independent. In each case we find that $ ( (\mathcal{F} \varphi_{1} )( 2 \pi \alpha)     )_{\alpha \in \mathbb{Z}}     $ is a vector that only contains zeros. So the result follows from Lemma \ref{lem_L2stab}.

\textit{Step 4.} Next we prove (iii) for $ m \in \mathbb{N}  $ odd with $ m \geq 3   $ and $p > 0$.  For the proof we use a generalization of the Strang-Fix conditions, see Theorem 4.2 in \cite{Ji1} and its proof. Of course the cardinal B-spline $N_{m}$ is compactly supported and refinable, see Lemma \ref{Bspl_elem}. Moreover for $m \geq 3$ it belongs to the Sobolev space $W^{1}_{1}(\mathbb{R})$. It is $(\mathcal{F}N_{m})(0) \not = 0$, see Lemma \ref{lem_qua_fourtrans}. Consequently Theorem 4.2 from \cite{Ji1} and the observation in its proof can be applied. With other words we have $ ( \mathcal{F} \varphi_{0} )^{(t)} (2 \pi k) = 0   $  for $ t \in \{ 0 , 1 \}   $ and $k \in \mathbb{Z} \setminus \{ 0 \}$. Therefore for $t = 0$ the sequence
\begin{equation}\label{strangfix1}
( \ldots , (\mathcal{F} \varphi_{0} )(- 4 \pi) , (\mathcal{F} \varphi_{0} )(- 2 \pi) , (\mathcal{F} \varphi_{0} )(0) , (\mathcal{F} \varphi_{0} )(2 \pi) , (\mathcal{F} \varphi_{0} )( 4 \pi) , \ldots )
\end{equation}
has the form $( \ldots , 0 , 0 , c_{0,m} , 0 , 0 , \ldots )$. For $t = 1$ we calculate
\begin{align*}
0 & = \Big ( \mathcal{F} \Big (   N_{m} \Big ( x + \lfloor \frac{m}{2} \rfloor  \Big )   \Big ) \Big )^{(1)}(2 \pi k) \\
& = - i \lceil \frac{m}{2} \rceil  \Big ( \mathcal{F} \Big (  \frac{x}{\lceil \frac{m}{2} \rceil }   N_{m} \Big ( x + \lfloor \frac{m}{2} \rfloor  \Big )   \Big ) \Big )(2 \pi k) \\
& = - i \lceil \frac{m}{2} \rceil   ( \mathcal{F} \varphi_{1}   )(2 \pi k) 
\end{align*}
with $ k \in \mathbb{Z} \setminus \{ 0 \}    $. Hence the sequence 
\begin{equation}\label{strangfix2}
( \ldots , (\mathcal{F} \varphi_{1} )(- 4 \pi) , (\mathcal{F} \varphi_{1} )(- 2 \pi) , (\mathcal{F} \varphi_{1} )(0) , (\mathcal{F} \varphi_{1} )(2 \pi) , (\mathcal{F} \varphi_{1} )( 4 \pi) , \ldots )
\end{equation}
has the form $( \ldots , 0 , 0 , c_{1,m} , 0 , 0 , \ldots )$. Combining \eqref{strangfix1} and \eqref{strangfix2} an application of Lemma \ref{lem_L2stab} yields that the integer translates of $ \Phi = (\varphi_{0}, \varphi_{1})^{T}$ are not $L_{2}-$stable. Therefore the claim follows.
\end{proof}

Now we summarize our previous findings by the following theorem. It is the main result of this section. 

\begin{Theorem}\label{thm_sec3summ1}
Let $ m \in \mathbb{N} $ and $ p \in \mathbb{N}_{0} $.  The sequence $  \{ V_{p,j}  \}_{j \in \mathbb{Z}}    $ of spaces from Definition \ref{def_Vpj} fulfills the properties (i)-(v) from Definition \ref{def_MRA}.
Let in addition either $ m= 1 $ or $  p=0   $. Then the sequence $  \{ V_{p,j}  \}_{j \in \mathbb{Z}}    $ is a multiresolution analysis according to Definition \ref{def_MRA} with a multiscaling function of multiplicity $p+1$.

\end{Theorem}

\begin{proof}
This result follows from Theorem \ref{thm_sec3summ}  and Lemma \ref{lem_L2stab_multiquark}.
\end{proof}

\section{The Perfect Reconstruction Property}\label{sec_perrec}

One of the goals of this paper is to identify suitable dual quarks and quarklets.  To this end in what follows we show that the quarks $ \Phi$ and the quarklets
\begin{equation}\label{quarklet_vec0}
\Psi := ( \psi_{0} , \psi_{1} , \ldots , \psi_{p} )^{T} 
\end{equation}
satisfy the so-called perfect reconstruction property. To see this at first let us recall some basic concepts concerning the theory of multiwavelets. The starting point of the theory of biorthogonal multiwavelets is the existence of a refinable function vector $   F = ( f_{1} , f_{2} , \ldots , f_{N} )^{T} \in ( L_{2}(\mathbb{R}))^{N}   $ with $  N \in \mathbb{N}  $. Recall that this vector is called \emph{refinable} if for $ - \infty < k_{0} \leq k_{1} < \infty   $ there exists a sequence of matrices $ \{ A_{k} \}_{k = k_{0}}^{k_{1}}  $ such that for all  $  x \in \mathbb{R}  $ we have
\begin{equation}\label{eq_ref_1}
\left ( \begin{array}{c} f_{1}(x) \\ \vdots  \\ f_{N}(x)  \end{array}   \right ) = \sum_{k = k_{0}}^{k_{1}} A_{k} \left ( \begin{array}{c} f_{1}(2x - k) \\ \vdots  \\ f_{N}(2x - k)  \end{array}   \right ) .
\end{equation}
By applying the Fourier transform we obtain
\begin{equation}\label{eq_matrix1}
\left ( \begin{array}{c} (\mathcal{F}f_{1})(\xi) \\ \vdots  \\ (\mathcal{F}f_{N})(\xi)  \end{array}   \right ) = \frac{1}{2} \sum_{k = k_{0}}^{k_{1}} A_{k} e^{- i k \frac{\xi}{2}} \left ( \begin{array}{c} (\mathcal{F}f_{1})(\frac{\xi}{2}) \\ \vdots  \\ (\mathcal{F}f_{N})(\frac{\xi}{2})  \end{array}   \right ) .
\end{equation}
We put $  z = e^{- i  \frac{\xi}{2}}  $ and define the \textit{symbol matrix} of $F$ via
\begin{equation}\label{eq_znot1}
\mathcal{A}(z) := \frac{1}{2} \sum_{k = k_{0}}^{k_{1}} A_{k} z^{k} .
\end{equation}
Then \eqref{eq_matrix1} can be rewritten as
\begin{equation}\label{eq_matrix2}
\left ( \begin{array}{c} (\mathcal{F}f_{1})(\xi) \\ \vdots  \\ (\mathcal{F}f_{N})(\xi)  \end{array}   \right ) = \mathcal{A}(z)   \left ( \begin{array}{c} (\mathcal{F}f_{1})(\frac{\xi}{2}) \\ \vdots  \\ (\mathcal{F}f_{N})(\frac{\xi}{2})  \end{array}   \right ) .
\end{equation}
A key issue of the theory of biorthogonal multiwavelets is the search for a second refinable function vector $   \tilde{F} = ( \tilde{f_{1}} , \tilde{f_{2}} , \ldots ,\tilde{f_{N}} )^{T} \in ( L_{2}(\mathbb{R}) )^{N}  $ which is biorthogonal (or dual) to $F$ in the following sense. 

\begin{Definition}\label{def_vec_biorto}
Let $   F = ( f_{1} , f_{2} , \ldots , f_{N} )^{T} \in ( L_{2}(\mathbb{R}))^{N}   $ and $   \tilde{F} = ( \tilde{f_{1}} , \tilde{f_{2}} , \ldots ,\tilde{f_{N}} )^{T} \in ( L_{2}(\mathbb{R}))^{N}   $ be two function vectors. Then they are called \emph{biorthogonal} if for all $  k \in \mathbb{Z}   $ 
\begin{equation}\label{bior_def_part11}
\Big ( 
\langle f_{i}( \cdot )  , \tilde{f_{j}}( \cdot - k )      \rangle_{L_{2}(\mathbb{R})}   \Big )_{i = 1, \ldots , N; j=1, \ldots ,N} = \delta_{0,k} Id_{N} .
\end{equation}
Here $ Id_{N}   $ is the identity matrix of size $  N \times N   $. 
\end{Definition}

For $ - \infty < k_{2} \leq k_{3} < \infty    $ the refinability of $ \tilde{F}  $ yields the existence of a sequence of matrices $ \{ \tilde{A_{k}} \}_{k = k_{2}}^{k_{3}}  $ such that for all  $  x \in \mathbb{R}  $ we have
\begin{equation}\label{eq_matrix3}
\left ( \begin{array}{c} \tilde{f_{1}}(x) \\ \vdots  \\ \tilde{f_{N}}(x)  \end{array}   \right ) = \sum_{k = k_{2}}^{k_{3}} \tilde{A_{k}} \left ( \begin{array}{c} \tilde{f_{1}}(2x - k) \\ \vdots  \\ \tilde{f_{N}}(2x - k)  \end{array}   \right ) .
\end{equation}
Again we can define the symbol matrix via
\begin{equation}\label{eq_znot2}
\tilde{\mathcal{A}}(z) := \frac{1}{2} \sum_{k = k_{2}}^{k_{3}} \tilde{A_{k}} z^{k} .
\end{equation}
Consequently \eqref{eq_matrix3} can be rewritten as
\begin{equation}\label{eq_matrix4}
\left ( \begin{array}{c} (\mathcal{F}\tilde{f_{1}})(\xi) \\ \vdots  \\ (\mathcal{F}\tilde{f_{N}})(\xi)  \end{array}   \right ) = \tilde{\mathcal{A}}(z)   \left ( \begin{array}{c} (\mathcal{F}\tilde{f_{1}})(\frac{\xi}{2}) \\ \vdots  \\ (\mathcal{F}\tilde{f_{N}})(\frac{\xi}{2})  \end{array}   \right ) .
\end{equation}
The biorthogonality of $F$ and $\tilde{F}$ yields an interesting connection between the symbols $ \mathcal{A}(z)   $ and $  \tilde{\mathcal{A}}(z)  $ which finally leads to the perfect reconstruction condition. It is given in the following well-known result, see Theorem 1 in \cite{YCW}. 

\begin{Lemma}\label{lem_perrec}
Let $   F = ( f_{1} , f_{2} , \ldots , f_{N} )^{T} \in ( L_{2}(\mathbb{R}) )^{N}   $ with $  N \in \mathbb{N}  $  be a function vector such that \eqref{eq_ref_1} is fulfilled. Let $   \tilde{F} = ( \tilde{f_{1}} , \tilde{f_{2}} , \ldots ,\tilde{f_{N}} )^{T} \in ( L_{2}(\mathbb{R}) )^{N}  $   be a second function vector such that \eqref{eq_matrix3} holds. For $  z = e^{- i  \frac{\xi}{2}}  $ we define the symbols $ \mathcal{A}(z)   $ and $  \tilde{\mathcal{A}}(z)  $ via \eqref{eq_znot1} and \eqref{eq_znot2}. If $ F $ and $ \tilde{F}  $ are biorthogonal, then we have
\begin{equation}\label{eq_pre_perfrec9}
\mathcal{A}(z) \overline{ \tilde{\mathcal{A}}(z) }^{T} + \mathcal{A}(- z)  \overline{ \tilde{\mathcal{A}}(- z) }^{T} = Id_{N} .
\end{equation}
\end{Lemma}

For a given function vector $F$ a popular method to find a biorthogonal vector $\tilde{F}$ is to look for a matrix $  \tilde{\mathcal{A}}(z)  $ such that \eqref{eq_pre_perfrec9} holds. For details concerning this topic we refer to \cite{Stre1} and \cite{StrStra}. Now let us turn to the quarks. We already know that the quark vector $  \Phi  $ is refinable, see Proposition \ref{prop_ref1}. Therefore the question arises whether there exists a dual vector $ \tilde{\Phi}   $ such that $ \Phi $ and $ \tilde{\Phi} $ are biorthogonal. However it turns out that in many cases there cannot exist a vector $ \tilde{\Phi} $ consisting of $ L_{2}(\mathbb{R})$ functions such that $ \Phi $ and $ \tilde{\Phi} $ are biorthogonal in the $L_{2}-$sense. So there is the following result. 

\begin{Theorem}\label{satz_biorto_not}
Let $m>1$ and $p>0$. Let  $ \Phi =  ( \varphi_{0} , \varphi_{1} , \ldots , \varphi_{p}    )^{T}    $ be the quark vector. Then there does not exist any function vector $ \tilde{ \Phi } =  ( \tilde{ \varphi_{0} } , \tilde { \varphi_{1} } , \ldots , \tilde {\varphi_{p}  }  )^{T} \in ( L_{2}(\mathbb{R}))^{p+1}    $ such that $ \Phi  $ and $ \tilde{ \Phi } $ are biorthogonal in the sense of Definition \ref{def_vec_biorto}.
\end{Theorem}

\begin{proof}
To prove this result we use Lemma 11.21 from \cite{Kein1}. It tells us that if $  \Phi, \tilde{\Phi}  \in ( L_{2}(\mathbb{R}))^{p+1}  $ are biorthogonal, they have stable shifts. However in Lemma \ref{lem_L2stab_multiquark} we learned that if $m$ and $p$ are as in Theorem \ref{satz_biorto_not}, the quark vector $ \Phi $ is not $L_{2}-$stable. 
\end{proof}

In what follows we have two goals. On the one hand we want to find out whether in the remaining case $m=1$ and $p \in \mathbb{N}_{0}$ which is not covered by Theorem \ref{satz_biorto_not} there exists a vector $ \tilde{\Phi} $ consisting of $ L_{2}(\mathbb{R})$ functions such that $ \Phi $ and $ \tilde{\Phi} $ are biorthogonal. On the other hand we intend to investigate whether for arbitrary $m$ and $p$ there exists a vector $\tilde{\Phi}$ that is dual to $\Phi$ in a more general way. We already mentioned that one possible way to answer such questions is to look for a matrix $  \tilde{\mathcal{A}}(z)  $ such that \eqref{eq_pre_perfrec9} is fulfilled. Before we start with this we have to incorporate the quarklets $ \Psi $ into our considerations. For that purpose we have to extend our definition of biorthogonality.

\begin{Definition}\label{def_biortomultiwavelet}
We call $   F = ( f_{1} , f_{2} , \ldots , f_{N} )^{T} \in ( L_{2}(\mathbb{R}) )^{N}   $ and $   \tilde{F} = ( \tilde{f}_{1} , \tilde{f}_{2} , \ldots ,\tilde{f}_{N} )^{T} \in ( L_{2}(\mathbb{R}) )^{N}   $ a pair of \emph{biorthogonal} function vectors if for all $  k \in \mathbb{Z}   $ the equation \eqref{bior_def_part11} holds. Moreover $   W = ( w_{1} , w_{2} , \ldots , w_{N} )^{T} \in ( L_{2}(\mathbb{R}) )^{N}   $ and $   \tilde{W} = ( \tilde{w}_{1} , \tilde{w}_{2} , \ldots ,\tilde{w}_{N} )^{T} \in ( L_{2}(\mathbb{R}) )^{N}   $ will be said to be a pair of \emph{biorthogonal} function vectors associated with $F$ and $\tilde{F}$, if for all $  k \in \mathbb{Z}   $ we have
\begin{align*}
\Big ( 
\langle f_{i}( \cdot )  , \tilde{w}_{j}( \cdot - k )      \rangle_{L_{2}(\mathbb{R})}  \Big )_{i=1 , \ldots , N; j = 1, \ldots , N} = 0_{N} 
\end{align*} 
and 
\begin{align*}
\Big ( 
\langle w_{i}( \cdot )  , \tilde{f}_{j}( \cdot - k )      \rangle_{L_{2}(\mathbb{R})}   \Big )_{i=1 , \ldots , N; j = 1, \ldots , N}  = 0_{N} 
\end{align*}
and
\begin{align*}
\Big ( 
\langle w_{i}( \cdot )  , \tilde{w}_{j}( \cdot - k )      \rangle_{L_{2}(\mathbb{R})}  \Big )_{i=1 , \ldots , N; j = 1, \ldots , N}  = \delta_{0,k} Id_{N} .
\end{align*}
Here $ 0_{N}   $ is the zero matrix.
\end{Definition}

Of course Definition \ref{def_biortomultiwavelet} is an extension of Definition \ref{def_vec_biorto}. So similar as in Lemma \ref{lem_perrec} biorthogonality yields an interesting property concerning the corresponding symbol matrices which can be found in the following lemma. 

\begin{Lemma}\label{lem_biortho_polyphasecond}
Let $   F = ( f_{1} , f_{2} , \ldots , f_{N} )^{T} \in ( L_{2}(\mathbb{R}) )^{N}   $ and $   \tilde{F} = ( \tilde{f_{1}} , \tilde{f_{2}} , \ldots ,\tilde{f}_{N} )^{T} \in ( L_{2}(\mathbb{R}) )^{N}  $ be a pair of compactly supported biorthogonal multiscaling functions that satisfy  \eqref{eq_ref_1} and \eqref{eq_matrix3}.
The matrices $  \mathcal{A}(z)  $ and $ \tilde{\mathcal{A}}(z)    $ are given by \eqref{eq_znot1} and \eqref{eq_znot2}. Moreover let $   W = ( w_{1} , w_{2} , \ldots , w_{N} )^{T} \in ( L_{2}(\mathbb{R}))^{N}   $ and $   \tilde{W} = ( \tilde{w}_{1} , \tilde{w}_{2} , \ldots ,\tilde{w}_{N} )^{T} \in ( L_{2}(\mathbb{R}) )^{N}  $  be a pair of biorthogonal multiwavelets associated with the multiscaling functions $F$ and $\tilde{F}$. We have
\begin{align*}
\left ( \begin{array}{c} w_{1}(x) \\ \vdots  \\ w_{N}(x)  \end{array}   \right ) = \sum_{k \in \mathbb{Z} } B_{k} \left ( \begin{array}{c} f_{1}(2x - k) \\ \vdots  \\ f_{N}(2x - k)  \end{array}   \right ) 
\end{align*}
and 
\begin{align*}
\left ( \begin{array}{c} \tilde{w}_{1}(x) \\ \vdots  \\ \tilde{w}_{N}(x)  \end{array}   \right ) = \sum_{k  \in \mathbb{Z} } \tilde{B}_{k} \left ( \begin{array}{c} \tilde{f_{1}}(2x - k) \\ \vdots  \\ \tilde{f}_{N}(2x - k)  \end{array}   \right ) .
\end{align*}
We put
\begin{equation}\label{eq_def_matB}
 \mathcal{B}(z) = \frac{1}{2} \sum_{k \in \mathbb{Z}} B_{k} z^{k} \qquad \mbox{and} \qquad \tilde{\mathcal{B}}(z) = \frac{1}{2} \sum_{k \in \mathbb{Z}} \tilde{B}_{k} z^{k} .  
\end{equation} 
Then for all $ z \in \mathbb{C}   $ with $ |z| = 1    $ we have 
\begin{equation}\label{eq_def_perfrec}
X(z) \overline{\tilde{X}(z)}^{T} := \left ( \begin{array}{cc} \mathcal{A}(z)   &   \mathcal{A}(-z)   \\
 \mathcal{B}(z)           &     \mathcal{B}(-z)      \\ \end{array}  \right )    \left ( \begin{array}{cc} \overline{ \tilde{\mathcal{A}}(z) }^{T}  &     \overline{ \tilde{\mathcal{B}}(z) }^{T} \\
\overline{ \tilde{\mathcal{A}}(-z) }^{T}          &   \overline{  \tilde{\mathcal{B}}(-z) }^{T}     \\ \end{array}  \right ) = Id_{2N} .
\end{equation}
\end{Lemma}

\begin{proof}
This result can be found in \cite{YCW}, see Theorem 1. Here we also can refer to (1.2) in \cite{GohJiaXia} or to formula (7.10) in \cite{Kein1}.
\end{proof}

The matrices $ X(z) $ and $  \overline{\tilde{X}(z)}^{T}   $   from Lemma \ref{lem_biortho_polyphasecond} are called \emph{modulation matrices}. The equation \eqref{eq_def_perfrec} describes the so-called \emph{perfect reconstruction property} which is essential in the theory of biorthogonal multiwavelets, see \cite{Kein1}, \cite{Stre1} and \cite{StrStra}. Therefore in what follows we investigate whether for the quarks and quarklets the condition of perfect reconstruction is fulfilled. Once we have found a matrix $ \tilde{X}(z)  $ for the quarks and quarklets it can be used to look for duals. Let us start with writing down the matrix $ X(z) $ for our case.

\begin{Lemma}\label{lem_qua_bio1}
Let $ m \in \mathbb{N}  $ and $ \tilde{m} \in \mathbb{N}    $ with $ m \leq \tilde{m}   $ and $ m + \tilde{m} \in 2 \mathbb{N}   $. Let $ p \in \mathbb{N}_{0}  $.   Let $   \Phi =  ( \varphi_{0} , \varphi_{1} , \ldots , \varphi_{p}    )^{T}        $ and $   \Psi = ( \psi_{0} , \psi_{1} , \ldots , \psi_{p} )^{T}      $. Then the matrix 
\begin{align*}
X(z) := \left ( \begin{array}{cc} \mathcal{A}(z)   &   \mathcal{A}(-z)   \\
 \mathcal{B}(z)           &     \mathcal{B}(-z)      \\ \end{array}  \right ) 
\end{align*}
from Lemma \ref{lem_biortho_polyphasecond} has the form  
\begin{align*}
X(z) = \left ( \begin{array}{cc}  \frac{1}{2} \sum_{k \in \mathbb{Z}} A_{k} z^{k}  &   \frac{1}{2} \sum_{k \in \mathbb{Z}} A_{k} (-z)^{k}   \\
 b(z) Id_{p+1}           &     b(-z) Id_{p+1}     \\ \end{array}  \right ) .
\end{align*}
Here the $A_{k}$ are the same as in Proposition \ref{prop_ref1}. It is $ b(z) := \frac{1}{2} \sum_{k \in \mathbb{Z}} b_{k} z^{k}   $, where the $b_{k}$ are the same as in Definition \ref{def_quarklet}.
\end{Lemma}

\begin{proof}
The shape of the upper part of the matrix follows from \eqref{eq_znot1} and Proposition \ref{prop_ref1}. For the lower part recall that $  \mathcal{B}(z)  $ is given by \eqref{eq_def_matB} with matrices $ \{ B_{k} \}_{k \in \mathbb{Z}}    $ defined via
\begin{align*}
\left ( \begin{array}{c} \psi_{0}(x) \\ \vdots  \\ \psi_{p}(x)  \end{array}   \right ) = \sum_{k \in \mathbb{Z} } B_{k} \left ( \begin{array}{c} \varphi_{0}(2x - k) \\ \vdots  \\ \varphi_{p}(2x - k)  \end{array}   \right ) .
\end{align*}
From Definition \ref{def_quarklet} we know
\begin{align*}
\psi_{p} = \sum_{k \in \mathbb{Z}} b_{k} \varphi_{p}(2 \cdot - k).
\end{align*}
Consequently we find $ B_{k} = b_{k} Id_{p+1}     $. 
\end{proof}

\begin{example}\label{ex_qua_polyphase}
Let $  m = \tilde{m} = 1    $ and $   \Phi =  ( \varphi_{0} , \varphi_{1}   )^{T}        $ and $   \Psi = ( \psi_{0} , \psi_{1}  )^{T}      $. Then we have
\begin{align*}
X(z) = \left ( \begin{array}{cccc}
\frac{1}{2} + \frac{1}{2} z  & 0  & \frac{1}{2} - \frac{1}{2} z  & 0   \\
\frac{1}{4} z  &  \frac{1}{4} + \frac{1}{4} z  & - \frac{1}{4} z  &  \frac{1}{4} - \frac{1}{4} z \\
\frac{1}{2} - \frac{1}{2} z  & 0  & \frac{1}{2} + \frac{1}{2} z  & 0 \\
 0  & \frac{1}{2} - \frac{1}{2} z  & 0  & \frac{1}{2} + \frac{1}{2} z  \\   
\end{array}  \right ) . 
\end{align*} 
\end{example}

It is already known that for the quarks and quarklets the modulation matrix $X(z)$ is invertible. There is the following result, see Proposition 4.14 in \cite{bib:sieber2020adaptive}. 

\begin{Lemma}\label{lem_qua_bio2}
Let $ m \in \mathbb{N}  $ and $ \tilde{m} \in \mathbb{N}    $ with $ m \leq \tilde{m}   $ and $ m + \tilde{m} \in 2 \mathbb{N}   $. Let $ p \in \mathbb{N}_{0}  $. Let $   \Phi =  ( \varphi_{0} , \varphi_{1} , \ldots , \varphi_{p}    )^{T}        $ and $   \Psi = ( \psi_{0} , \psi_{1} , \ldots , \psi_{p} )^{T}      $. For $  z \in \mathbb{C}   $ with $ |z|=1   $ we define $T(z) := \mathcal{A}(z) b(-z) - b(z) \mathcal{A}(-z)$. Then $T(z)^{-1}$ exists and the matrix $  X(z)  $ from Lemma \ref{lem_qua_bio1} is invertible with
\begin{align*}
X(z)^{-1} =  \left ( \begin{array}{cc} b(-z) T(z)^{-1}  &  - T(z)^{-1}    \mathcal{A}(-z)  \\
   - b(z) T(z)^{-1}      &  T(z)^{-1} \mathcal{A}(z)      \\ \end{array}  \right ) .
\end{align*}
\end{Lemma}

Here we also can refer to \cite{bib:DRS19}, where a slightly different notation was used.  

\begin{example}\label{ex_qua_polyphase2}
Let $  m = \tilde{m} = 1    $ and $   \Phi =  ( \varphi_{0} , \varphi_{1}   )^{T}        $ and $   \Psi = ( \psi_{0} , \psi_{1}  )^{T}      $. Then the matrix $  X(z)  $ is given in Example \ref{ex_qua_polyphase}.  For the inverse matrix $  X(z)^{-1}  $ we obtain
\begin{align*}
X(z)^{-1} = \frac{1}{2} \left ( \begin{array}{cccc}
z^{-1} + 1  & 0  & - z^{-1} + 1  & 0   \\
- \frac{1}{2} z^{-1} - \frac{1}{2}  &  2 z^{-1} + 2  &  \frac{1}{2} z^{-1}  + \frac{1}{2} &  - z^{-1} + 1 \\
- z^{-1} + 1  & 0  & z^{-1} + 1  & 0 \\
 \frac{1}{2} z^{-1} - \frac{1}{2}  & - 2 z^{-1} + 2  &  - \frac{1}{2} z^{-1} + \frac{1}{2}  & z^{-1} + 1  \\   
\end{array}  \right ) . 
\end{align*} 
\end{example}

In what follows we will prove that the matrix $ X(z)^{-1} $ can be interpreted as dual modulation matrix $ \tilde{X}(z) $, see \eqref{eq_def_perfrec}. For that purpose we need some preparations. 

\begin{Lemma}\label{lem_qua_bio3}
Let $ m \in \mathbb{N}  $ and $ \tilde{m} \in \mathbb{N}    $ with $ m \leq \tilde{m}   $ and $ m + \tilde{m} \in 2 \mathbb{N}   $. Let $ p \in \mathbb{N}_{0}  $. Let $   \Phi =  ( \varphi_{0} , \varphi_{1} , \ldots , \varphi_{p}    )^{T}        $ and $   \Psi = ( \psi_{0} , \psi_{1} , \ldots , \psi_{p} )^{T}      $. Let $  z \in \mathbb{C}   $ with $ |z|=1   $. We use the notation from Lemma \ref{lem_qua_bio2}. Then the following assertions are true. 
\begin{itemize}
\item[(i)] The matrix $ T(z) $ is of lower triangular shape. It only consists of Laurent polynomials and fulfills $  T(-z) = - T(z)   $. 

\item[(ii)] The matrix $ T(z)^{-1} $ is of lower triangular shape. It only consists of Laurent polynomials and fulfills $  T(-z)^{-1} = - T(z)^{-1}      $.

\item[(iii)] It is $ b(z) T(-z)^{-1} = - b(z) T(z)^{-1} $ and $ - T(-z)^{-1}    \mathcal{A}(z) =   T(z)^{-1} \mathcal{A}(z)   $. 

\end{itemize}
\end{Lemma}

\begin{proof}
\textit{Step 1.} At first we prove (i). The lower triangular shape of $ T(z) $ has already been observed in \cite{bib:sieber2020adaptive}, see the proof of Proposition 4.14. It is clear that both $  \mathcal{A}(z)  $ and $ b(z) $ only consist of Laurent polynomials. Concerning the last statement we observe
\begin{align*}
- T(z) & = - [ \mathcal{A}(z) b(-z) - b(z) \mathcal{A}(-z)   ] =   \mathcal{A}(-z) b(z) - b(-z) \mathcal{A}(z)  = T(-z) .
\end{align*}
Therefore the claim follows. 

\textit{Step 2.} Now we prove (ii). The fact that $ T(z)^{-1} $ is of lower triangular shape follows from the lower triangular shape of $ T(z)$. It already has been observed in the Appendix of \cite{bib:DRS19} that $ T(z)^{-1} $ only consists of Laurent polynomials. Moreover we can use (i) to find
\begin{align*}
T(-z)^{-1} = [ - T(z) ]^{-1} = (-1)^{-1} T(z)^{-1} = - T(z)^{-1} .
\end{align*}
So the proof of (ii) is complete.

\textit{Step 3.} Now we prove (iii). For that purpose at first we observe $   b(z) T(-z)^{-1} = - b(z) T(z)^{-1} $. This directly follows from (ii). So we have proved the first part of the assertion. To continue we have to show $  - T(-z)^{-1}    \mathcal{A}(z) =   T(z)^{-1} \mathcal{A}(z)   $. Let us write $ T(-z)^{-1}    \mathcal{A}(z) = D(z)   $. Then the entry $ d_{i,j}(z)$ of $ D(z) $ has the form
\begin{align*}
d_{i,j}(z) = \sum_{k = 1}^{p+1} (T(-z)^{-1} )_{i,k}  ( \mathcal{A}(z) )_{k,j} = - \sum_{k = 1}^{p+1} (T(z)^{-1} )_{i,k}  ( \mathcal{A}(z) )_{k,j} .
\end{align*}
Here we also used (ii). Consequently we can conclude
\begin{align*}
- T(-z)^{-1}    \mathcal{A}(z) =   T(z)^{-1} \mathcal{A}(z) .
\end{align*} 
The proof is complete. 
\end{proof}

Now we are well-prepared to show that quarks and quarklets fulfill the perfect reconstruction condition. 

\begin{Theorem}\label{Satz_qua_perref1}
Let $ m \in \mathbb{N}  $ and $ \tilde{m} \in \mathbb{N}    $ with $ m \leq \tilde{m}   $ and $ m + \tilde{m} \in 2 \mathbb{N}   $. Let $ p \in \mathbb{N}_{0}  $. Let $   \Phi =  ( \varphi_{0} , \varphi_{1} , \ldots , \varphi_{p}    )^{T}        $ and $   \Psi = ( \psi_{0} , \psi_{1} , \ldots , \psi_{p} )^{T}      $. Let $  z \in \mathbb{C}   $ with $ |z|=1   $. Let $X(z)$ be as in Lemma \ref{lem_qua_bio1}. Then $ X(z)   $ is invertible and for the inverse we observe
\begin{align*}
X(z)^{-1} = \overline{ \tilde{X}(z) }^{T} = \left ( \begin{array}{cc} \overline{ \tilde{\mathcal{A}}(z) }^{T}  &     \overline{ \tilde{\mathcal{B}}(z) }^{T} \\
\overline{ \tilde{\mathcal{A}}(-z) }^{T}          &   \overline{  \tilde{\mathcal{B}}(-z) }^{T}     \\ \end{array}  \right ) .
\end{align*}
Here we put
\begin{equation}\label{eq_mathA_mathB1}
\overline{\tilde{\mathcal{A}}(z) }^{T} =  b(-z) T(z)^{-1} \qquad \mbox{and} \qquad \overline{ \tilde{\mathcal{B}}(z) }^{T} = - T(z)^{-1}    \mathcal{A}(-z) .
\end{equation}
In other words the perfect reconstruction condition \eqref{eq_def_perfrec} is fulfilled. Moreover there exist finitely supported matrix sequences  $ \{ \tilde{A}_{k} \}_{k \in \mathbb{Z}}  $ and $ \{ \tilde{B}_{k} \}_{k \in \mathbb{Z}}  $ such that 
\begin{equation}\label{eq_mathA_mathB_sum1}
\tilde{\mathcal{A}}(z) = \frac{1}{2} \sum_{k \in \mathbb{Z}} \tilde{A}_{k} z^{k}  \qquad \mbox{and}  \qquad \tilde{\mathcal{B}}(z) = \frac{1}{2} \sum_{k \in \mathbb{Z}} \tilde{B}_{k} z^{k} .  
\end{equation} 
\end{Theorem}

\begin{proof}
For the proof we can collect the things we did before. That $ X(z) $ is invertible follows from Lemma \ref{lem_qua_bio2}. The special structure of $   X(z)^{-1} $ was shown in the Lemmas \ref{lem_qua_bio2} and \ref{lem_qua_bio3}. It remains to prove that the matrices $ \tilde{\mathcal{A}}(z)  $ and $  \tilde{\mathcal{B}}(z)  $ can be decomposed into sequences $  \{ \tilde{A}_{k} \}_{k \in \mathbb{Z}}    $ and $   \{ \tilde{B}_{k} \}_{k \in \mathbb{Z}}    $. For that purpose we have to examine \eqref{eq_mathA_mathB1}. From Lemma \ref{lem_qua_bio3} we know that $ T(z)^{-1}  $ only consists of Laurent polynomials. Moreover it is $ b(z) = \frac{1}{2} \sum_{k \in \mathbb{Z}} b_{k} z^{k}   $ and $\mathcal{A}(z)$ is defined as in Lemma \ref{lem_qua_bio1}. It is known that sums and products of Laurent polynomials again are Laurent polynomials. Hence $ \overline{\tilde{\mathcal{A}}(z) }  $ and $    \overline{\tilde{\mathcal{B}}(z) } $ only consist of Laurent polynomials.  Let $  z \in \mathbb{C}   $ with $ |z|=1   $. Then for each Laurent polynomial we observe
\begin{align*}
\overline{p(z)} = \overline{ \sum_{k = n_{1}}^{n_{2}} a_{k} z^{k} } = \sum_{k = n_{1}}^{n_{2}} \overline{ a_{k} } \overline{ z^{k} } =  \sum_{k = n_{1}}^{n_{2}} \overline{ a_{k} } ( \overline{ z } )^{k} = \sum_{k = n_{1}}^{n_{2}} \overline{ a_{k} }  z^{-k} .
\end{align*}
Here we used that since $ |z|=1   $ we have $\overline{z} =  z^{-1} $. Hence also $ \tilde{\mathcal{A}}(z)  $ and $    \tilde{\mathcal{B}}(z) $ only consist of Laurent polynomials. Thus there exist matrix sequences  $ \{ \tilde{A}_{k} \}_{k \in \mathbb{Z}}  $ and $ \{ \tilde{B}_{k} \}_{k \in \mathbb{Z}}  $ such that \eqref{eq_mathA_mathB_sum1} holds. 
\end{proof}

\begin{example}\label{ex_qua_polyphase4}
Let $  m = \tilde{m} = 1    $ and $   \Phi =  ( \varphi_{0} , \varphi_{1}   )^{T}        $ and $   \Psi = ( \psi_{0} , \psi_{1}  )^{T}      $. Then for $  \tilde{\mathcal{A}}(z)    $ and $  \tilde{\mathcal{B}}(z)   $ we obtain
\begin{align*}
 \tilde{\mathcal{A}}(z)    =  \left ( \begin{array}{cc} \frac{1}{2} z + \frac{1}{2}  &  - \frac{1}{4} z - \frac{1}{4}  \\
    0     &  z  + 1  \\ \end{array}  \right )  \qquad \mbox{and} \qquad  \tilde{\mathcal{B}}(z)   =  \left ( \begin{array}{cc} - \frac{1}{2} z + \frac{1}{2}  &   \frac{1}{4} z + \frac{1}{4}   \\
    0    & - \frac{1}{2} z  + \frac{1}{2}   \\ \end{array}  \right ) .
\end{align*}
Consequently the matrix sequences  $ \{ \tilde{A}_{k} \}_{k \in \mathbb{Z}}  $ and $ \{ \tilde{B}_{k} \}_{k \in \mathbb{Z}}  $ are given by
\begin{align*}
 \tilde{A}_{0}   =  \left ( \begin{array}{cc} 1  &  - \frac{1}{2}    \\
    0     & 2  \\ \end{array}  \right )  \qquad \mbox{and} \qquad   \tilde{A}_{1}   =  \left ( \begin{array}{cc} 1  & -  \frac{1}{2}  \\
    0    & 2   \\ \end{array}  \right ) 
\end{align*}
and
\begin{align*}
 \tilde{B}_{0}   =  \left ( \begin{array}{cc} 1  &  \frac{1}{2}    \\
    0     & 1  \\ \end{array}  \right )  \qquad \mbox{and} \qquad   \tilde{B}_{1}   =  \left ( \begin{array}{cc} -1  &   \frac{1}{2}  \\
    0    & -1   \\ \end{array}  \right ). 
\end{align*}
For $ k \in \mathbb{Z} \setminus \{ 0 , 1 \}    $ we have $ \tilde{A}_{k}       = \tilde{B}_{k } = 0_{2}  $. Here $ 0_{2}  $ is the zero matrix.

\end{example}

\section{Generalized dual Quarklets}\label{sec_dualqua}

In Section \ref{sec_perrec} we have seen that quarks and quarklets fulfill the perfect reconstruction property. Moreover we obtained a dual modulation matrix $  \tilde{X}(z)   $. This knowledge can now be used to establish results concerning the existence and the properties of generalized dual quarks and generalized dual quarklets. For that purpose we apply some results from the theory of biorthogonal multiwavelets which we recall in the following. Here the subsequent definition will be important. 

\begin{Definition}\label{def_CondE}
We say that a square matrix satisfies \emph{Condition E}, if it has a simple eigenvalue of 1, and all other eigenvalues are smaller than 1 in absolute value. 
\end{Definition}

Condition E is very important within the theory of biorthogonal multiwavelets, see for example \cite{PlStr} or \cite{Kein1}. It can be used to formulate the following result which is already known. 

\begin{Lemma}\label{lem_conv_matrix}
Let $  \tilde{\mathcal{A}}(z)  $ be defined via \eqref{eq_znot2} with $ z = e^{- i  \frac{\xi}{2}} $. Then the following assertions are true. 

\begin{itemize}
\item[(i)] The matrix refinement equation \eqref{eq_matrix3} has a compactly supported distributional solution vector $   \tilde{F} = ( \tilde{f_{1}} , \tilde{f_{2}} , \ldots ,\tilde{f}_{N} )^{T} \in ( \mathcal{S}'(\mathbb{R}) )^{N}   $ if and only if $  \tilde{\mathcal{A}}(1)   $ has an eigenvalue of the form $   2^{n}$ with $ n \in \mathbb{N}_{0}    $. 

\item[(ii)] Let $1$ be the only eigenvalue of $   \tilde{\mathcal{A}}(1)   $ of the form $   2^{n}$ with $ n \in \mathbb{N}_{0}    $. Then the matrix refinement equation \eqref{eq_matrix3} has a unique compactly supported distributional solution (up to a constant)  $   \tilde{F} = ( \tilde{f_{1}} , \tilde{f_{2}} , \ldots ,\tilde{f}_{N} )^{T}    $ with $ ( (\mathcal{F} \tilde{f_{1}})(0) , (\mathcal{F}\tilde{f_{2}})(0)  , \ldots , (\mathcal{F}\tilde{f}_{N})(0)  )^{T} = v      $ if and only if $1$ is a simple eigenvalue of $   \tilde{\mathcal{A}}(1)   $ and $  v =     \tilde{\mathcal{A}}(1) v        $. In particular we have
\begin{equation}\label{eq_iter_dualsc2}
\left ( \begin{array}{c} (\mathcal{F}\tilde{f_{1}})(\xi) \\ \vdots  \\ (\mathcal{F}\tilde{f}_{N})(\xi)  \end{array}   \right ) = \prod_{j = 1}^{\infty}  \tilde{\mathcal{A}}( e^{- i 2^{-j}  \xi } )   \left ( \begin{array}{c} (\mathcal{F}\tilde{f_{1}})( 0 ) \\ \vdots  \\ (\mathcal{F}\tilde{f}_{N})( 0 )  \end{array}   \right )   .
\end{equation}

\item[(iii)] Assume that the function vector $   \tilde{F} = ( \tilde{f_{1}} , \tilde{f_{2}} , \ldots ,\tilde{f}_{N} )^{T} \in (L_{2}(\mathbb{R}))^{N}    $ is compactly supported, refinable and $L_{2}-$stable. Then $   \tilde{\mathcal{A}}(1)   $ satisfies Condition E.
\end{itemize}
\end{Lemma}

\begin{proof}
All results can be found in \cite{PlStr}, see Theorems 1-3. Let us also refer to \cite{JianShen} and \cite{Kein1}, see Chapter 11.
\end{proof}

Now we can use Lemma \ref{lem_conv_matrix} to prove the existence of generalized dual quarks. More precisely there is the following result. 

\begin{Lemma}\label{lem_qua_dual1}
Let $ m \in \mathbb{N}  $ and $ \tilde{m} \in \mathbb{N}    $ with $ m \leq \tilde{m}   $ and $ m + \tilde{m} \in 2 \mathbb{N}   $. Let $ p \in \mathbb{N}_{0}  $ and $   \Phi =  ( \varphi_{0} , \varphi_{1} , \ldots , \varphi_{p}    )^{T}        $ and $   \Psi = ( \psi_{0} , \psi_{1} , \ldots , \psi_{p} )^{T}      $. Let $ \{ \tilde{A}_{k} \}_{k \in \mathbb{Z}} $ be the matrix sequence from Theorem \ref{Satz_qua_perref1} with \eqref{eq_znot2}. Then the matrix refinement equation 
\begin{equation}\label{eq_notL2prof1}
\left ( \begin{array}{c} \tilde{\varphi}_{0} \\ \vdots  \\ \tilde{\varphi}_{p}  \end{array}   \right ) = \sum_{k \in \mathbb{Z}} \tilde{A}_{k} \left ( \begin{array}{c} \tilde{\varphi}_{0}(2 \cdot - k) \\ \vdots  \\ \tilde{\varphi}_{p}(2 \cdot - k)  \end{array}   \right ) 
\end{equation}
has a compactly supported distributional solution vector $   \tilde{\Phi} = ( \tilde{\varphi}_{0} , \tilde{\varphi}_{1} , \ldots ,\tilde{\varphi}_{p} )^{T} \in ( \mathcal{S}'(\mathbb{R}) )^{p+1}   $.
\end{Lemma}

\begin{proof}
To prove this result we use part (i) of Lemma \ref{lem_conv_matrix}. That means we have to compute the eigenvalues of $ \tilde{\mathcal{A}}(1)   $. Recall that we have
\begin{align*}
\overline{\tilde{\mathcal{A}}(z) }^{T} =  b(-z) T(z)^{-1} .
\end{align*}
We know that $ T(z)^{-1}  $ has lower triangular shape, see Lemma \ref{lem_qua_bio3}. Moreover $b(z)$ is a Laurent polynomial. Consequently to calculate the eigenvalues of $ \tilde{\mathcal{A}}(1)  $ we have to know the diagonal entries of $  T(z)^{-1}  $. We already have seen that also $ T(z)  $ is of lower triangular shape with diagonal entries
\begin{align*}
( T(z) )_{q,q} = 2^{-q+1}  z ,
\end{align*}
see the proof of Proposition 4.14 in \cite{bib:sieber2020adaptive}. Here we have $ q \in \{ 1 , 2 , \ldots , p+1 \}    $. Hence $ T(z)^{-1}   $ has diagonal entries of the form 
\begin{align*}
( T(z)^{-1} )_{q,q} = 2^{q-1} \frac{1}{  z } .
\end{align*}
Now let us turn to the Laurent polynomial $ b(z) $. From the theory of biorthogonal wavelets (see for example Section 2.1 in \cite{bib:sieber2020adaptive} or Remark 3.6 in \cite{bib:DRS19}) we know that $ b(-z) $ can be written in the form 
\begin{align*}
b(-z)  =  z \overline{\tilde{a}(z)} =  z \overline{\frac{1}{2} \sum_{k \in \mathbb{Z}} \tilde{a}_{k} z^{k}} = z \frac{1}{2} \sum_{k \in \mathbb{Z}} \tilde{a}_{k} \overline{z}^{k} .
\end{align*}
Here the numbers $ \tilde{a}_{k} \in \mathbb{R}  $ are well-known and explicitly given in \cite{bib:sieber2020adaptive}, see Theorem 2.7. Let us also refer to Section 6.A in \cite{bib:CDF92}. Recall that we are interested in $ \tilde{\mathcal{A}}(1)    $. 
We observe 
\begin{align*}
b(-1)  =   \overline{\tilde{a}(1)}  =  \frac{1}{2} \sum_{k \in \mathbb{Z}} \tilde{a}_{k}  .
\end{align*}
Fortunately this sum is well-known from the theory of CDF-wavelets, see Section 6.A on page 541 in \cite{bib:CDF92}. Here with a slightly different notation we learn $ b(-1) =  \overline{\tilde{a}(1)} = 1   $. Hence we get
\begin{align*}
\tilde{\mathcal{A}}(1)  = \overline{  T(1)^{-1} }^{T} .
\end{align*}
Recall that $  T(1)^{-1}  $ is a lower triangular matrix with diagonal entries $( T(1)^{-1} )_{q,q} = 2^{q-1}  $. Consequently for the diagonal entries of $\tilde{\mathcal{A}}(1) $ we observe
\begin{align*}
( \tilde{\mathcal{A}}(1) )_{q,q} = ( \overline{  T(1)^{-1} }^{T} )_{q,q} = ( T(1)^{-1} )_{q,q} = 2^{q-1}  .
\end{align*}
Since $  T(1)^{-1}   $ is a lower triangular matrix $ \tilde{\mathcal{A}}(1)  $ is an upper triangular matrix. Consequently $ \tilde{\mathcal{A}}(1)  $ has $ p + 1   $ eigenvalues of the form
\begin{align*}
\lambda_{1} = 1, \qquad \lambda_{2} = 2, \qquad \ldots \qquad \lambda_{p+1} = 2^{p} .
\end{align*}
Therefore for all $ p \in \mathbb{N}_{0}   $ there is an eigenvalue of the form $ 2^{n} $ with $ n \in \mathbb{N}_{0}   $. Hence we can apply part (i) of Lemma \ref{lem_conv_matrix}. So the claim follows.  
\end{proof}

Next we want to study the properties of the generalized dual quarks $  \tilde{\Phi}  $ we found in Lemma \ref{lem_qua_dual1}.  

\begin{Lemma}\label{lem_qua_dual2}
Let $ m \in \mathbb{N}  $ and $ \tilde{m} \in \mathbb{N}    $ with $ m \leq \tilde{m}   $ and $ m + \tilde{m} \in 2 \mathbb{N}   $. Let $ p \in \mathbb{N}  $ with $ p \not = 0   $ and $   \Phi =  ( \varphi_{0} , \varphi_{1} , \ldots , \varphi_{p}    )^{T}        $ and $   \Psi = ( \psi_{0} , \psi_{1} , \ldots , \psi_{p} )^{T}      $. Let $ \{ \tilde{A}_{k} \}_{k \in \mathbb{Z}} $ be the matrix sequence from Theorem \ref{Satz_qua_perref1} with \eqref{eq_znot2}. Then the compactly supported distributional solution vector $   \tilde{\Phi} = ( \tilde{\varphi}_{0} , \tilde{\varphi}_{1} , \ldots ,\tilde{\varphi}_{p} )^{T} \in ( \mathcal{S}'(\mathbb{R})   )^{p+1} $ of the matrix refinement equation \eqref{eq_notL2prof1} is not $L_{2}-$stable.
\end{Lemma}

\begin{proof}
We prove this result via contradiction. Assume that the solution vector $   \tilde{\Phi} = ( \tilde{\varphi}_{0} , \tilde{\varphi}_{1} , \ldots ,\tilde{\varphi}_{p} )^{T} \in ( \mathcal{S}'(\mathbb{R}) )^{p+1}   $ also belongs to $ (  L_{2}(\mathbb{R})   )^{p+1}   $ and is $  L_{2}-  $stable. Then we can apply part (iii) of Lemma \ref{lem_conv_matrix}. Therefore we find that $   \tilde{\mathcal{A}}(1)   $ has a simple eigenvalue $1$ and the moduli of all its other eigenvalues are less then $1$. From the proof of Lemma \ref{lem_qua_dual1} we know that $ \tilde{\mathcal{A}}(1)  $ has $ p + 1   $ eigenvalues of the form
\begin{align*}
\lambda_{1} = 1, \qquad \lambda_{2} = 2, \qquad \ldots \qquad \lambda_{p+1} = 2^{p} .
\end{align*}
Since $ p > 0   $ consequently $2$ is an eigenvalue of $ \tilde{\mathcal{A}}(1)    $. But this is a contradiction.  
\end{proof}

Unfortunately it turns out that in most of the cases the generalized dual quarks which are collected in the vector $ \tilde{\Phi}   $ are not biorthogonal to $\Phi$ in the sense of Definition \ref{def_vec_biorto}. 

\begin{Theorem}\label{satz_dual_notL2}
Let $ m \in \mathbb{N}  $ and $ \tilde{m} \in \mathbb{N}    $ with $ m \leq \tilde{m}   $ and $ m + \tilde{m} \in 2 \mathbb{N}   $. Let $ p \in \mathbb{N}  $ with $ p \not = 0   $ and $   \Phi =  ( \varphi_{0} , \varphi_{1} , \ldots , \varphi_{p}    )^{T}        $ and $   \Psi = ( \psi_{0} , \psi_{1} , \ldots , \psi_{p} )^{T}      $. 
Then the compactly supported distributional solution vector $   \tilde{\Phi} = ( \tilde{\varphi}_{0} , \tilde{\varphi}_{1} , \ldots ,\tilde{\varphi}_{p} )^{T} \in ( \mathcal{S}'(\mathbb{R}) )^{p+1}   $ of the matrix refinement equation \eqref{eq_notL2prof1} is not biorthogonal to $ \Phi  $ in the sense of Definition \ref{def_vec_biorto}.
\end{Theorem}

\begin{proof}
We prove this result by contradiction. Assume that the solution vector $ \tilde{ \Phi } =  ( \tilde{ \varphi}_{0}  , \tilde { \varphi}_{1}  , \ldots , \tilde {\varphi}_{p}    )^{T}    $ belongs to $ ( L_{2}(\mathbb{R}) )^{p+1}  $ and $ \Phi  $ and $ \tilde{ \Phi } $ are biorthogonal. Now because of $ \Phi =  ( \varphi_{0} , \varphi_{1} , \ldots , \varphi_{p}    )^{T}   \in ( L_{2}(\mathbb{R}) )^{p+1}   $ we can apply Lemma 11.21 from \cite{Kein1}. It tells us that $ \tilde{ \Phi } $ is $L_{2}-$stable. But this is a contradiction to Lemma \ref{lem_qua_dual2}.
\end{proof}

Nevertheless it is possible to describe the generalized dual quarks explicitly via an infinite product. 

\begin{Lemma}\label{lem_qua_dualconst1}
Let $ m \in \mathbb{N}  $ and $ \tilde{m} \in \mathbb{N}    $ with $ m \leq \tilde{m}   $ and $ m + \tilde{m} \in 2 \mathbb{N}   $. Let $ p \in \mathbb{N}_{0}  $ and $   \Phi =  ( \varphi_{0} , \varphi_{1} , \ldots , \varphi_{p}    )^{T}        $ and $   \Psi = ( \psi_{0} , \psi_{1} , \ldots , \psi_{p} )^{T}      $. Let $ \{ \tilde{A}_{k} \}_{k \in \mathbb{Z}} $ be the matrix sequence from Theorem \ref{Satz_qua_perref1} with \eqref{eq_znot2}. Let $v$ be a right eigenvector of $ 2^{-p} \tilde{\mathcal{A}}(1)    $ corresponding to the eigenvalue $1$. By $ D^{p}    $ we denote the p-th distributional derivative operator. Then the matrix refinement equation \eqref{eq_notL2prof1} has a compactly supported distributional solution vector $   \tilde{\Phi} = ( \tilde{\varphi}_{0} , \tilde{\varphi}_{1} , \ldots ,\tilde{\varphi}_{p} )^{T} \in ( \mathcal{S}'(\mathbb{R}) )^{p+1}   $ which is given by
\begin{align*}
&\left ( \begin{array}{c} \tilde{\varphi}_{0} \\ \vdots  \\ \tilde{\varphi}_{p}  \end{array}   \right )  = D^{p} \mathcal{F}^{-1} \Big [ \prod_{j = 1}^{\infty} 2^{-p} \tilde{\mathcal{A}}( e^{- i 2^{-j}  \xi } ) v \Big ] .
\end{align*}
\end{Lemma}

\begin{proof}
For the proof we use some ideas from \cite{JianShen}, see Theorem 2.5. At first let us recall that $ \tilde{\mathcal{A}}(1)  $ has $ p + 1   $ eigenvalues of the form
\begin{align*}
\lambda_{1} = 1, \qquad \lambda_{2} = 2, \qquad \ldots \qquad \lambda_{p+1} = 2^{p} .
\end{align*}
Now we put $\tilde{\mathcal{A}}_{0}(z) := 2^{-p} \tilde{\mathcal{A}}(z) $. Since $ \tilde{\mathcal{A}}(1)  $ is an upper triangular matrix also $ \tilde{\mathcal{A}}_{0}(1)   $ is an upper triangular matrix. For the diagonal entries we find
\begin{align*}
( \tilde{\mathcal{A}}_{0}(1) )_{q,q}  = 2^{q-1-p}  .
\end{align*}
Here we have $ q \in \{ 1 , 2 , \ldots , p+1 \}    $.  Consequently $ \tilde{\mathcal{A}}_{0}(1)  $ has $ p + 1   $ eigenvalues of the form
\begin{align*}
\mu_{1} = 2^{-p}, \qquad \mu_{2} = 2^{-p+1}, \qquad \ldots \qquad \mu_{p+1} = 1 .
\end{align*}
Hence we can apply part (ii) of Lemma \ref{lem_conv_matrix}. Let $v$ be a right eigenvector of $ \tilde{\mathcal{A}}_{0}(1)    $ corresponding to the eigenvalue $1$. Then there exists a unique compactly supported distribution vector (up to a constant)  $   \tilde{F} = ( \tilde{f}_{0} , \tilde{f}_{1} , \ldots ,\tilde{f}_{p} )^{T}    $ with $ ( (\mathcal{F} \tilde{f}_{0})(0) , (\mathcal{F}\tilde{f}_{1})(0)  , \ldots , (\mathcal{F}\tilde{f}_{p})(0)  )^{T} = v      $ such that 
\begin{align*}
\left ( \begin{array}{c} (\mathcal{F}\tilde{f}_{0})(\xi) \\ \vdots  \\ (\mathcal{F}\tilde{f}_{p})(\xi)  \end{array}   \right ) = \tilde{\mathcal{A}}_{0}( e^{- i  \frac{\xi}{2}} )   \left ( \begin{array}{c} (\mathcal{F}\tilde{f}_{0})(\frac{\xi}{2}) \\ \vdots  \\ (\mathcal{F}\tilde{f}_{p})(\frac{\xi}{2})  \end{array}   \right ) .
\end{align*}
Moreover we have
\begin{align*}
\left ( \begin{array}{c} (\mathcal{F}\tilde{f_{0}})(\xi) \\ \vdots  \\ (\mathcal{F}\tilde{f_{p}})(\xi)  \end{array}   \right ) = \prod_{j = 1}^{\infty}  \tilde{\mathcal{A}}_{0}( e^{- i 2^{-j}  \xi } )   \left ( \begin{array}{c} (\mathcal{F}\tilde{f_{0}})( 0 ) \\ \vdots  \\ (\mathcal{F}\tilde{f_{p}})( 0 )  \end{array}   \right ) = \prod_{j = 1}^{\infty}  \tilde{\mathcal{A}}_{0}( e^{- i 2^{-j}  \xi } ) v   
\end{align*}
with pointwise convergence. With other words we have 
\begin{align*}
\left ( \begin{array}{c} \tilde{f}_{0} \\ \vdots  \\ \tilde{f}_{p} \end{array}   \right ) = \mathcal{F}^{-1} \Big [ \prod_{j = 1}^{\infty}  \tilde{\mathcal{A}}_{0}( e^{- i 2^{-j}  \xi } ) v \Big ] .
\end{align*}
Since the distributions $ \tilde{f}_{q-1} \in \mathcal{S}'(\mathbb{R})   $ are compactly supported, the Paley-Wiener-Schwartz Theorem yields that the $ \mathcal{F} \tilde{f}_{q-1} $ are entire analytic functions of at most polynomial growth. Now let $ D^{p}    $ be the p-th distributional derivative operator. Then we put 
\begin{align*}
\tilde{\Phi} = ( \tilde{\varphi}_{0} , \tilde{\varphi}_{1} , \ldots ,\tilde{\varphi}_{p} )^{T} := D^{p} \tilde{F} := ( D^{p} \tilde{f_{0}} , D^{p} \tilde{f_{1}} , \ldots , D^{p} \tilde{f_{p}} )^{T} . 
\end{align*}
Then of course $ \tilde{\Phi} \in ( \mathcal{S}'(\mathbb{R}) )^{p+1}   $ consists of compactly supported distributions. We observe
\begin{align*}
\left ( \begin{array}{c} (\mathcal{F} \tilde{\varphi}_{0} )(\xi) \\ \vdots  \\ (\mathcal{F} \tilde{\varphi}_{p}  )(\xi)  \end{array}   \right ) = C_{1} \left ( \begin{array}{c} \xi^{p} (\mathcal{F} \tilde{f_{0}} )(\xi) \\ \vdots  \\ \xi^{p} (\mathcal{F} \tilde{f_{p}}  )(\xi)  \end{array}   \right ) \quad \mbox{and} \quad  \frac{ 2^{p} \xi^{-p} }{C_{1}} \left ( \begin{array}{c} (\mathcal{F} \tilde{\varphi}_{0} )(\frac{\xi}{2}) \\ \vdots  \\ (\mathcal{F} \tilde{\varphi}_{p}  )(\frac{\xi}{2})  \end{array}   \right ) =  \left ( \begin{array}{c}  (\mathcal{F} \tilde{f_{0}} )(\frac{\xi}{2}) \\ \vdots  \\  (\mathcal{F} \tilde{f_{p}}  )(\frac{\xi}{2})  \end{array}   \right ) .
\end{align*}
Moreover we have
\begin{align*}
\left ( \begin{array}{c} (\mathcal{F} \tilde{\varphi}_{0} )(\xi) \\ \vdots  \\ (\mathcal{F} \tilde{\varphi}_{p}  )(\xi)  \end{array}   \right )  = C_{1} \xi^{p} \tilde{\mathcal{A}}_{0}( e^{- i  \frac{\xi}{2}} )   \left ( \begin{array}{c} (\mathcal{F}\tilde{f_{0}})(\frac{\xi}{2}) \\ \vdots  \\ (\mathcal{F}\tilde{f_{p}})(\frac{\xi}{2})  \end{array}   \right ) =  2^{p} \tilde{\mathcal{A}}_{0}( e^{- i  \frac{\xi}{2}} )     \left ( \begin{array}{c} (\mathcal{F} \tilde{\varphi}_{0} )(\frac{\xi}{2}) \\ \vdots  \\ (\mathcal{F} \tilde{\varphi}_{p}  )(\frac{\xi}{2})  \end{array}   \right )  .
\end{align*}
Now recall that we have $ \tilde{\mathcal{A}}_{0}(z) = 2^{-p} \tilde{\mathcal{A}}(z)   $. Consequently we obtain 
\begin{align*}
\left ( \begin{array}{c} (\mathcal{F} \tilde{\varphi}_{0} )(\xi) \\ \vdots  \\ (\mathcal{F} \tilde{\varphi}_{p}  )(\xi)  \end{array}   \right )   =   \tilde{\mathcal{A}}( e^{- i  \frac{\xi}{2}} )     \left ( \begin{array}{c} (\mathcal{F} \tilde{\varphi}_{0} )(\frac{\xi}{2}) \\ \vdots  \\ (\mathcal{F} \tilde{\varphi}_{p}  )(\frac{\xi}{2})  \end{array}   \right )  .
\end{align*}
Therefore we get
\begin{align*}
\left ( \begin{array}{c} \tilde{\varphi}_{0} \\ \vdots  \\ \tilde{\varphi}_{p}  \end{array}   \right ) = \sum_{k \in \mathbb{Z}} \tilde{A}_{k} \left ( \begin{array}{c} \tilde{\varphi}_{0}(2 \cdot - k) \\ \vdots  \\ \tilde{\varphi}_{p}(2 \cdot - k)  \end{array}   \right ). 
\end{align*}
Moreover our construction yields that $ \tilde{\Phi} $ can be written as 
\begin{align*}
&\left ( \begin{array}{c} \tilde{\varphi}_{0} \\ \vdots  \\ \tilde{\varphi}_{p}  \end{array}   \right ) = D^{p} \left ( \begin{array}{c} \tilde{f}_{0} \\ \vdots  \\ \tilde{f}_{p}  \end{array}   \right )  = D^{p} \mathcal{F}^{-1} \Big [ \prod_{j = 1}^{\infty} 2^{-p} \tilde{\mathcal{A}}( e^{- i 2^{-j}  \xi } ) v \Big ] .
\end{align*}
So the proof is complete. 
\end{proof}

Using the results we obtained up to now it is also possible to define generalized dual quarklets. For instance this can be done in the following way. 

\begin{Definition}\label{def_dual_quarkl}
Let $ m \in \mathbb{N}  $ and $ \tilde{m} \in \mathbb{N}    $ with $ m \leq \tilde{m}   $ and $ m + \tilde{m} \in 2 \mathbb{N}   $. Let $ p \in \mathbb{N}_{0}  $ and $   \Phi =  ( \varphi_{0} , \varphi_{1} , \ldots , \varphi_{p}    )^{T}        $ and $   \Psi = ( \psi_{0} , \psi_{1} , \ldots , \psi_{p} )^{T}      $. Moreover   $ \{ \tilde{A}_{k} \}_{k \in \mathbb{Z}} $ is the matrix sequence from Theorem \ref{Satz_qua_perref1}. Let $   \tilde{\Phi} = ( \tilde{\varphi}_{0} , \tilde{\varphi}_{1} , \ldots ,\tilde{\varphi}_{p} )^{T} \in ( \mathcal{S}'(\mathbb{R}) )^{p+1}  $ be the compactly supported distributional solution vector of  the matrix refinement equation \eqref{eq_notL2prof1} given in Lemma \ref{lem_qua_dualconst1}. Furthermore $ \{ \tilde{B}_{k} \}_{k \in \mathbb{Z}}  $ is the matrix sequence from Theorem \ref{Satz_qua_perref1}. Then we define the \emph{generalized dual quarklets} $  \tilde{\Psi} = ( \tilde{\psi}_{0} , \tilde{\psi}_{1} , \ldots , \tilde{\psi}_{p} )^{T}      $ via 
\begin{equation}\label{eq_dual_qu}
\left ( \begin{array}{c} \tilde{\psi}_{0} \\ \vdots  \\ \tilde{\psi}_{p}  \end{array}   \right ) = \sum_{k \in \mathbb{Z}} \tilde{B}_{k} \left ( \begin{array}{c} \tilde{\varphi}_{0}(2 \cdot - k) \\ \vdots  \\ \tilde{\varphi}_{p}(2 \cdot - k)  \end{array}   \right ) .
\end{equation}
\end{Definition}

It is not difficult to see that the properties of the generalized dual quarklets are similar to those of the generalized dual quarks. More precisely we know the following. 

\begin{Lemma}\label{lem_dualqua_prop1}
Let $ m \in \mathbb{N}  $ and $ \tilde{m} \in \mathbb{N}    $ with $ m \leq \tilde{m}   $ and $ m + \tilde{m} \in 2 \mathbb{N}   $. Let $ p \in \mathbb{N}_{0}  $. Moreover $  \tilde{\Psi} = ( \tilde{\psi}_{0} , \tilde{\psi}_{1} , \ldots , \tilde{\psi}_{p} )^{T}      $ are the dual quarklets from Definition \ref{def_dual_quarkl}. Then all $  \tilde{\psi}_{q}  $ with $  q \in \{ 0 , 1 , \ldots , p \}   $ are compactly supported distributions from $ \mathcal{S}'(\mathbb{R})    $. 
\end{Lemma}

\begin{proof}
From Definition \ref{def_dual_quarkl} we know that the generalized dual quarklets $  \tilde{\Psi} = ( \tilde{\psi}_{0} , \tilde{\psi}_{1} , \ldots , \tilde{\psi}_{p} )^{T}      $ are given by \eqref{eq_dual_qu}. In Lemma \ref{lem_qua_dualconst1} we learned that $   \tilde{\Phi} = ( \tilde{\varphi}_{0} , \tilde{\varphi}_{1} , \ldots ,\tilde{\varphi}_{p} )^{T} \in ( \mathcal{S}'(\mathbb{R}) )^{p+1}   $ consists of compactly supported distributions. The matrix sequence $ \{ \tilde{B}_{k} \}_{k \in \mathbb{Z}}  $ is given in Theorem \ref{Satz_qua_perref1}. Here we observed that 
\begin{align*}
\tilde{\mathcal{B}}(z) = \frac{1}{2} \sum_{k \in \mathbb{Z}} \tilde{B}_{k} z^{k} \qquad \qquad \mbox{and} \qquad \qquad \overline{ \tilde{\mathcal{B}}(z) }^{T} = - T(z)^{-1}    \mathcal{A}(-z) .
\end{align*}
From Proposition \ref{prop_ref1} we know that $ \mathcal{A}(z)   $ consists of a finite number of matrices $A_{k}$. In Lemma \ref{lem_qua_bio3} we learned that $  T(z)^{-1}   $ only consists of Laurent polynomials with finitely many non-vanishing coefficients. Consequently the same holds true for $ T(z)^{-1}    \mathcal{A}(-z)  $. Therefore the claim follows.  
\end{proof}
 
\begin{Remark}
Recall that for $p \not = 0$ in Theorem \ref{satz_dual_notL2} we learned that $  \Phi$ and $ \tilde{\Phi}  $ are not biorthogonal in the usual $ L_{2}-$sense, see Definition \ref{def_vec_biorto}. However there might be a chance that $ \Phi  $ and $ \tilde{\Phi}  $ are dual in a more general way at least in some special cases. So we know that our quarks $ \Phi  $ belong to the Sobolev space $ ( H^{s}(\mathbb{R}) )^{p+1}  $ for $ s < m - \frac{1}{2}  $, see Lemma 3 in Chapter  2.3.1 in \cite{RS}. Consequently it arises the question whether also $ \tilde{\Phi} $ belongs to a Sobolev space $ ( H^{t}(\mathbb{R}) )^{p+1}  $ for some $ t \in \mathbb{R}  $, possibly with $ t < 0  $. One possible way to deal with this problem is to use the theory developed in \cite{Han}. More precisely in \cite{Han} we can apply Theorem 6.3.3 on page 517 in combination with (5.6.44) and Theorem 5.8.4 on page 461. Proceeding that way for example in the special case $  m = \tilde{m} = 1    $ with $ p = 1 $ and $   \Phi =  ( \varphi_{0} , \varphi_{1}   )^{T}        $ we obtain $ \tilde{\Phi}  \in (H^{-1}(\mathbb{R}))^{2} $. However maybe this result is not sharp, see \cite{Han}. To obtain optimal results for arbitrary parameters $m, \tilde{m}$ and $p$ much more effort is necessary. This will be the subject of future research. Finally one can ask whether $ \Phi $ and $ \tilde{\Phi}  $ are dual in the sense of dual Sobolev spaces, see Chapter 4.6 on page 327 in \cite{Han} for a definition. However duality can not hold in the case $ m > 1   $ and $ p > 0   $ because of the lack of stability, see the proof of Theorem 6.4.5 in \cite{Han}. So biorthogonality would imply $L_{2}-$stability of $\Phi$, what contradicts Lemma \ref{lem_L2stab_multiquark}.
\end{Remark}

\section{Quarklet Multiscale Decomposition}\label{sec_decompoqua}

In this section we want to find out whether quarks and quarklets provide a multiscale decomposition of $  L_{2}(\mathbb{R})  $. The starting point will be the spaces $ V_{p,j}   $, see Definition \ref{def_Vpj}. We want to know whether there exist complement spaces $  W_{p,j}  $ generated by the quarklets such that $  V_{p, j+1} = V_{p,j} \oplus W_{p,j}    $. Here in general $ \oplus$ can not be expected to be an orthogonal sum. But it might be a non-orthogonal algebraic direct sum. One first step towards the answer is the following lemma. Here we see that each vector of fine quarks has a decomposition into quarklets and quarks on a coarse level.  

\begin{Lemma}\label{quarklet_ref1}
Let $ m \in \mathbb{N}  $ and $ \tilde{m} \in \mathbb{N}    $ with $ m \leq \tilde{m}   $ and $ m + \tilde{m} \in 2 \mathbb{N}   $. Let $ p \in \mathbb{N}_{0}  $. Put $ \Phi =  ( \varphi_{0} , \varphi_{1} , \ldots , \varphi_{p}    )^{T}   $ and $ \Psi = ( \psi_{0} , \psi_{1} , \ldots , \psi_{p} )^{T}   $. Let $  \varrho \in \{ 0 , 1\}  $. Then there exist sequences of $ (p+1) \times (p+1) $ matrices $ \{ C_{\varrho + 2k} \}_{k \in \mathbb{Z}}   $ and $ \{ D_{\varrho + 2k} \}_{k \in \mathbb{Z}}   $ such that
\begin{align*}
\Phi ( 2 \cdot - \varrho ) = \sum_{k \in \mathbb{Z}} C_{\varrho + 2k} \Phi ( \cdot - k ) + \sum_{k \in \mathbb{Z}}  D_{\varrho + 2k} \Psi ( \cdot - k ) .
\end{align*}
Moreover the length of the reconstruction sequences  $ \{ C_{\varrho + 2k} \}_{k \in \mathbb{Z}}   $ and  $ \{ D_{\varrho + 2k} \}_{k \in \mathbb{Z}}   $ scales linearly with $p$.

\end{Lemma}

This result can be found in \cite{bib:sieber2020adaptive}, see Theorem 4.15 and Theorem 4.16. A similar observation also was made in \cite{bib:DRS19}. In order to prove the desired decomposition of $L_{2}(\mathbb{R})$ it is advantageous to introduce the concept of sub-symbol and polyphase matrices.   

\begin{Definition}\label{def_subsym_11}
Let $ m \in \mathbb{N}  $ and $ \tilde{m} \in \mathbb{N}    $ with $ m \leq \tilde{m}   $ and $ m + \tilde{m} \in 2 \mathbb{N}   $. Let $ p \in \mathbb{N}_{0}  $. Let $   \Phi =  ( \varphi_{0} , \varphi_{1} , \ldots , \varphi_{p}    )^{T}        $ and $   \Psi = ( \psi_{0} , \psi_{1} , \ldots , \psi_{p} )^{T}      $. The matrices $ \mathcal{A}(z)  $ and $ \mathcal{B}(z) $ as well as the matrix sequences $ \{ A_{k} \}_{k}     $ and $ \{B_{k} \}_{k}    $ are as in Lemma \ref{lem_qua_bio1}.

\begin{itemize}
\item[(i)] Then for $ \varrho \in \{ 0 , 1 \}    $ we define the \emph{sub-symbol matrices} $  \mathcal{A}_{\varrho}  (z^2) $ and $  \mathcal{B}_{\varrho}  (z^2) $  by
\begin{align*}
\left ( \mathcal{A}_{\varrho}  (z^2) \right )_{q,l} := \sum_{k \in \mathbb{Z}}  ( A_{2 k + \varrho} )_{q,l} z^{2 k} 
\end{align*}
and
\begin{align*}
\left ( \mathcal{B}_{\varrho}  (z^2) \right )_{q,l} := \sum_{k \in \mathbb{Z}}  ( B_{2 k + \varrho} )_{q,l} z^{2 k} 
\end{align*}
with $ q , l \in \{ 1 , \ldots , p + 1 \}      $.

\item[(ii)] We define the \emph{polyphase matrix} $ P(z) $ via
\begin{align*}
P(z) := \left ( \begin{array}{cc} \mathcal{A}_{0}  (z^2)  &  \mathcal{A}_{1}  (z^2)  \\
     \mathcal{B}_{0}  (z^2)    &  \mathcal{B}_{1}  (z^2) \\ \end{array}  \right )  .
\end{align*}
\end{itemize}

\end{Definition}

Definition \ref{def_subsym_11} was inspired by Proposition 3.1 in \cite{bib:DRS19}. The polyphase matrix plays a crucial role in the theory of biorthogonal multiwavelets, see Chapter 7.7 in \cite{Kein1} and  \cite{GohYap}. It can be written in terms of the modulation matrix, see Proposition 3.3 in \cite{bib:DRS19} for the quarklet case.

\begin{Lemma}\label{lem_polyphase_invert}
Let $ m \in \mathbb{N}  $ and $ \tilde{m} \in \mathbb{N}    $ with $ m \leq \tilde{m}   $ and $ m + \tilde{m} \in 2 \mathbb{N}   $. Let $ p \in \mathbb{N}_{0}  $. Let $   \Phi =  ( \varphi_{0} , \varphi_{1} , \ldots , \varphi_{p}    )^{T}        $ and $   \Psi = ( \psi_{0} , \psi_{1} , \ldots , \psi_{p} )^{T}      $. Then the polyphase matrix $ P(z)  $ is invertible for all $ z \in \mathbb{C}   $ with $ |z| = 1    $. In other words it has full rank of $  2p + 2  $.
\end{Lemma}

\begin{proof}
For the proof at first we observe that we can write
\begin{align*}
P(z) = X(z) \mathcal{E}(z)
\end{align*}
with $ X(z) $ as in Lemma \ref{lem_qua_bio1} and 
\begin{align*}
\mathcal{E}(z) := \left ( \begin{array}{cc} Id_{p+1}  &  \frac{1}{z} Id_{p+1} \\
     Id_{p+1}   &  - \frac{1}{z} Id_{p+1} \\ \end{array}  \right )  .
\end{align*}
This is proved in \cite{bib:DRS19}, see Proposition 3.3. There we also find that $  \mathcal{E}(z)  $ is invertible for all $ z \in \mathbb{C}   $ with $ |z| = 1    $. Moreover we know that $ X(z) $ is invertible, see Lemma \ref{lem_qua_bio2}. Consequently we find that also the polyphase matrix $ P(z) $ is invertible. 
\end{proof}

In what follows the spaces $  W_{p,j}  $ generated by the quarklets will be important. 

\begin{Definition}\label{def_Wpj1}
Let $ j \in \mathbb{Z}  $ and $  p \in \mathbb{N}_{0}  $. Then we define the spaces
\begin{align*}
W_{p,j} :=  \overline{ \spa \{ 2^{ \frac{j}{2} } \psi_{q} (2^{j} \cdot - k ) : 0 \leq q \leq p , k \in \mathbb{Z}   \} } \subset L_{2}(\mathbb{R}).
\end{align*}
Here the functions $ \psi_{q}  $ are the quarklets from Definition \ref{def_quarklet}.
\end{Definition}

The spaces $  W_{p,j}  $ can be used to obtain the desired decomposition of $L_{2}(\mathbb{R})$. One important intermediate step toward this result can be found in the following lemma. 

\begin{Lemma}\label{lem_oplus1}
Let $ m \in \mathbb{N}  $ and $ \tilde{m} \in \mathbb{N}    $ with $ m \leq \tilde{m}   $ and $ m + \tilde{m} \in 2 \mathbb{N}   $. Let $   p \in \mathbb{N}_{0}  $ and $j \in \mathbb{Z}$.

\begin{itemize}
\item[(i)]
Then we have $ V_{p, j+1} = V_{p,j} + W_{p,j} $.   
\item[(ii)]
Let in addition either $ m =  1 $ or $  p=0   $. Then we also have
\begin{align*}
V_{p, j+1} = V_{p,j} \oplus W_{p,j}  .   
\end{align*}
Here $  \oplus  $ denotes the (non-orthogonal) algebraic direct sum.
\end{itemize}
\end{Lemma}

\begin{proof}
For the proof we only consider the case $ j = 0  $. Then the general case follows by standard arguments.

\textit{Step 1.} At first we prove $  V_{p, 1} \subset V_{p,0} + W_{p,0}   $. For that purpose let $  f \in  V_{p, 1}  $. Then $f$ can be approximated in $  L_{2}(\mathbb{R}) $ via functions of the form
\begin{align*}
g = \sum_{q = 0}^{p} \sum_{l \in \mathbb{Z}} c_{l,q} \varphi_{q} (2 \cdot - l) 
\end{align*}
with appropriate $  \{ c_{l,q} \}_{l \in \mathbb{Z}, q \in \{ 0,1, \ldots , p \} } \subset \mathbb{C}  $. From Lemma \ref{quarklet_ref1} we know that for $ \varrho \in \{ 0 , 1 \}  $ and $  q \in \{   0 , 1 , \ldots , p \}   $ we can find sequences  $  \{ e_{k,r , \varrho} \}_{k \in \mathbb{Z}, r \in \{ 0,1, \ldots , p \} , \varrho \in \{ 0 ,1 \} } \subset \mathbb{C}  $ and $  \{ d_{k,r , \varrho} \}_{k \in \mathbb{Z}, r \in \{ 0,1, \ldots , p \} , \varrho \in \{ 0 ,1 \} } \subset \mathbb{C}  $  such that 
\begin{align*}
\varphi_{q} (2 \cdot - \varrho ) = \sum_{k \in \mathbb{Z}} \sum_{r = 0}^{p} e_{k,r,\varrho} \varphi_{r} ( \cdot - k ) +  \sum_{k \in \mathbb{Z}} \sum_{r = 0}^{p} d_{k,r, \varrho} \psi_{r} ( \cdot - k ) .
\end{align*}
For $  l \in \mathbb{Z} $ similar representations can be found for $ \varphi_{q} (2 \cdot - l)   $. Consequently we can write
\begin{align*}
g & = \sum_{q = 0}^{p} \sum_{l \in \mathbb{Z}} c_{l,q} \varphi_{q} (2 \cdot - l)  \\
& = \sum_{q = 0}^{p} \sum_{l \in \mathbb{Z}} c_{l,q} \Big (   \sum_{k \in \mathbb{Z}} \sum_{r = 0}^{p} e_{l,q,k,r} \varphi_{r} ( \cdot - k ) +  \sum_{k \in \mathbb{Z}} \sum_{r = 0}^{p} d_{l,q,k,r} \psi_{r} ( \cdot - k ) \Big )   \\
& = \sum_{r = 0}^{p} \sum_{k \in \mathbb{Z}}   \alpha_{k,r} \varphi_{r} ( \cdot - k ) +  \sum_{r = 0}^{p} \sum_{k \in \mathbb{Z}}    \beta_{k,r} \psi_{r} ( \cdot - k )    
\end{align*}
for appropriate $    \{ e_{l,q,k,r} \}_{l,q,k,r } \subset \mathbb{C}    $, $    \{ d_{l,q,k,r} \}_{l,q,k,r } \subset \mathbb{C}    $, $    \{ \alpha_{k,r} \}_{k,r } \subset \mathbb{C}    $ and $    \{ \beta_{k,r} \}_{k,r } \subset \mathbb{C}    $.   Therefore we conclude $  f \in V_{p,0} + W_{p,0}    $.

\textit{Step 2.} Now we prove $ V_{p,0} + W_{p,0} \subset   V_{p, 1}   $. Let $  f \in  V_{p,0} + W_{p,0}  $. Then $f$ can be approximated in $  L_{2}(\mathbb{R}) $ via functions of the form
\begin{align*}
g = \sum_{q = 0}^{p} \sum_{l \in \mathbb{Z}} d_{l,q} \varphi_{q} ( \cdot - l) + \sum_{q = 0}^{p} \sum_{l \in \mathbb{Z}} c_{l,q} \psi_{q} ( \cdot - l)  
\end{align*}
with appropriate $  \{ d_{l,q} \}_{l \in \mathbb{Z}, q \in \{ 0,1, \ldots , p \} } \subset \mathbb{C}  $ and  $  \{ c_{l,q} \}_{l \in \mathbb{Z}, q \in \{ 0,1, \ldots , p \} } \subset \mathbb{C}  $. We use Proposition \ref{prop_ref1} to find
\begin{align*}
g = \sum_{q = 0}^{p} \sum_{l \in \mathbb{Z}} a_{l,q} \varphi_{q} ( 2 \cdot - l) + \sum_{q = 0}^{p} \sum_{l \in \mathbb{Z}} c_{l,q} \psi_{q} ( \cdot - l)  
\end{align*}
with $  \{ a_{l,q} \}_{l \in \mathbb{Z}, q \in \{ 0,1, \ldots , p \} } \subset \mathbb{C}  $.  Then Definition \ref{def_quarklet} yields
\begin{align*}
g = \sum_{q = 0}^{p} \sum_{l \in \mathbb{Z}} a_{l,q} \varphi_{q} ( 2 \cdot -  l) + \sum_{q = 0}^{p} \sum_{l \in \mathbb{Z}} b_{l,q} \varphi_{q} ( 2 \cdot - l)  
\end{align*}
with $  \{ b_{l,q} \}_{l \in \mathbb{Z}, q \in \{ 0,1, \ldots , p \} } \subset \mathbb{C}  $.   Consequently we find $   f \in   V_{p, 1}  $. 

\textit{Step 3.} Now we prove (ii). For that purpose let in addition either $ m= 1 $ or $  p=0   $. Recall that for $p=0$ the desired result is already well-known, see Section 6.A in \cite{bib:CDF92}. So it remains to deal with the case $ m = 1 $. Here at first from Lemma \ref{lem_L2stab_multiquark} we learn that the integer shifts of $ \Phi =  ( \varphi_{0} , \varphi_{1} , \ldots , \varphi_{p}    )^{T}    $ are $L_{2}-$stable. Consequently we can use similar arguments as described in Section 1 in \cite{GohYap}, see also Section 4 in \cite{LeeTanTang}. Since the integer shifts of $ \Phi =  ( \varphi_{0} , \varphi_{1} , \ldots , \varphi_{p}    )^{T}    $ are $L_{2}-$stable and because of Definition \ref{def_Vpj} our spaces $ V_{p,0} $ fit into the setting described at the beginning of Section 4 in \cite{LeeTanTang}. Therefore we can use the arguments from \cite{LeeTanTang}, see formula (4.3), to find that the integer shifts of 
\begin{equation}\label{eq_stable_newpro}
\varphi_{0}(2x) , \varphi_{1}(2x) , \ldots , \varphi_{p}(2x), \varphi_{0}(2x-1) , \varphi_{1}(2x-1) , \ldots , \varphi_{p}(2x-1)
\end{equation}
are a Riesz basis for $ V_{p,1}   $. Those are $ 2p + 2 $ functions. Next from Definition \ref{def_quarklet} and Proposition \ref{prop_ref1} we conclude that
\begin{align*}
\{ \varphi_{0} , \varphi_{1} , \ldots , \varphi_{p}    , \psi_{0} , \psi_{1} , \ldots , \psi_{p} \} \subset V_{p,1} .
\end{align*}
We observe that the integer shifts of the $ 2p + 2 $ functions 
\begin{equation}\label{eq_proof_polyph_nonmo}
 \varphi_{0} , \varphi_{1} , \ldots , \varphi_{p}    , \psi_{0} , \psi_{1} , \ldots , \psi_{p} 
\end{equation}
are $L_{2}-$stable. To see this we can use Theorem 4.3 in \cite{JiMi}. It can be applied since the functions in \eqref{eq_stable_newpro} have stable integer translates. Recall that all involved functions belong to $  \mathcal{L}_{2}(\mathbb{R})     $, see the remark before Theorem 2.1 in \cite{JiMi} for a definition and explanations.  Moreover that the corresponding polyphase matrix $P(z)$ from Definition \ref{def_subsym_11} is nonsingular for all $z \in \mathbb{C}$ with $|z|=1$ has already been observed in Lemma \ref{lem_polyphase_invert}. Consequently the $L_{2}-$stability of \eqref{eq_proof_polyph_nonmo} follows. Let us remark that this result also can be obtained by using Theorem 5.1 from \cite{JiMi} and Lemma \ref{lem_L2stab}. However proceeding that way one has to carry out some rather technical computations. In a next step Theorem 3.1 from \cite{LeeTanTang} yields
\begin{align*}
\overline { \spa \{ \varphi_{0}( \cdot - k ) ,  \ldots , \varphi_{p}( \cdot - k )    , \psi_{0}( \cdot - k )  , \ldots , \psi_{p}( \cdot - k )  :  k \in \mathbb{Z} \} }  = V_{p,1} .
\end{align*}
Consequently we are in the setting described in \cite{GohYap}, see page 143. So like there we can conclude that $ V_{p,1} $ is the algebraic direct sum of $ V_{p,0} $ and $  W_{p,0}  $. This follows from Definition \ref{def_Wpj1} which is similar to that from \cite{GohYap}.
\end{proof}

Now we are well-prepared to prove the desired decomposition of $  L_{2}(\mathbb{R})   $ using the spaces $  W_{p,j}   $.

\begin{Theorem}\label{satz_L2zerl1}
Let $ m \in \mathbb{N}  $ and $ \tilde{m} \in \mathbb{N}    $ with $ m \leq \tilde{m}   $ and $ m + \tilde{m} \in 2 \mathbb{N}   $. Let $   p \in \mathbb{N}_{0}  $.
\begin{itemize}
\item[(i)] We have
\begin{align*}
L_{2}(\mathbb{R} ) & = \overline{ \ldots + W_{p,-2} + W_{p,-1} + W_{p,0} + W_{p,1} + W_{p,2} + \ldots } \\
& = \overline{ V_{p,0} + W_{p,0} + W_{p,1} + W_{p,2} + \ldots } \ .
\end{align*}

\item[(ii)] Let in addition either $ m = 1 $ or $  p=0   $. Then we have
\begin{align*}
L_{2}(\mathbb{R} ) = \overline{ \bigoplus_{j \in \mathbb{Z}} W_{p,j} } = \overline{ V_{p,0} \oplus \bigoplus_{j = 0}^{\infty} W_{p,j} } .
\end{align*}
Here $  \oplus  $ denotes the (non-orthogonal) algebraic direct sum.
\end{itemize}
\end{Theorem}

\begin{proof}
We only proof (ii). The proof of (i) can be done with similar methods using Lemma \ref{lem_oplus1}.

\textit{Step 1.} At first we prove that for all $ J \in \mathbb{N}   $ we have 
\begin{equation}\label{eq_prov_L2dec_new1}
\overline{V_{p,0} \oplus \bigoplus_{j = 0}^{J} W_{p,j} } \subset L_{2}(\mathbb{R}) .
\end{equation}
To see this we use Lemma \ref{lem_oplus1} several times. Then we get 
\begin{align*}
\overline{V_{p,0} \oplus \bigoplus_{j = 0}^{J} W_{p,j}} = \overline{ V_{p,1} \oplus \bigoplus_{j = 1}^{J} W_{p,j} } = \ldots = V_{p,J+1} \subset L_{2}(\mathbb{R}) 
\end{align*}
due to Definition \ref{def_Vpj}. This shows \eqref{eq_prov_L2dec_new1}.

\textit{Step 2.} Next we prove
\begin{equation}\label{eq_prov_L2dec_new2}
L_{2}(\mathbb{R}) \subset \overline{ V_{p,0} \oplus \bigoplus_{j = 0}^{\infty} W_{p,j} }  .
\end{equation}
For that purpose we apply Theorem \ref{thm_sec3summ} and Lemma \ref{lem_oplus1} to find
\begin{align*}
L_{2}(\mathbb{R}) = \overline{ \bigcup_{j \in \mathbb{Z}} V_{p,j} } =  \overline{ \bigcup_{j \geq 1} V_{p,j} } =  \overline{ \bigcup_{j \geq 1} V_{p,j-1} \oplus W_{p,j-1} }. 
\end{align*}
Iterating this argument we get 
\begin{align*}
L_{2}(\mathbb{R})   =  \overline{ \bigcup_{j \geq 1}   \bigoplus_{n=0}^{j-1} W_{p,n} \oplus V_{p,0} }  \subseteq      \overline{ \bigoplus_{n=0}^{\infty} W_{p,n} \oplus V_{p,0} } .   
\end{align*}
So the proof of \eqref{eq_prov_L2dec_new2} is complete.

\textit{Step 3.} It remains to show that we have 
\begin{equation}\label{eq_prov_L2dec_new3}
 \bigoplus_{j \in \mathbb{Z}} W_{p,j}  =  V_{p,0} \oplus \bigoplus_{j = 0}^{\infty} W_{p,j}  .
\end{equation}
To see this Lemma \ref{lem_oplus1} yields that for all $  J \in \mathbb{Z} \setminus \mathbb{N}_{0}$ we have 
\begin{align*}
 V_{p,0} \oplus \bigoplus_{j = 0}^{\infty} W_{p,j} = V_{p,J} \oplus \bigoplus_{j = J}^{\infty} W_{p,j} .    
\end{align*}
Sending $J$ to $- \infty $ and using Theorem \ref{thm_sec3summ}, see $(i)$ and $(iii)$ in Definition \ref{def_MRA}, finally \eqref{eq_prov_L2dec_new3} follows. This also completes the proof of Theorem \ref{satz_L2zerl1}.
\end{proof}

One very important application of quarks and quarklets is to use them for the numerical approximation of functions. For that purpose it seems to be advantageous if the functions that generate $  W_{p,j}  $ are orthogonal. Therefore in what follows we increase the orthogonality of those spaces at least for $m=1$. To this end we can use some results from \cite{Kaz1}. Here the concept of minimality (sometimes also called g-minimality) plays an important role. For that let us refer to page 2 in \cite{Kaz1}, see also Definition 2.4 in \cite{TaKaz}.

\begin{Definition}\label{def_gmini}
Let $  F = ( f_{1}, \ldots , f_{N}  )^{T} \in ( L_{2}(\mathbb{R}))^{N}  $ with $   N \in \mathbb{N}  $ be a finite vector of nontrivial functions. Then we use the following definitions.  
\begin{itemize}
\item[(i)] We put $S(F) :=   \overline{ \spa \{ f_{l} ( \cdot - k ) :  k \in \mathbb{Z}  , f_{l} \in F  \} } \subset L_{2}(\mathbb{R}) $.
\item[(ii)] For any $  f_{l} \in F  $ with $   1 \leq l \leq N  $ we put $   F^{(l)}  := F \setminus \{ f_{l}   \} $.
\item[(iii)] We say that $F$ is \emph{minimal} (or \emph{g-minimal}) if for any $  1 \leq l \leq N  $ we have $  f_{l} \not \in S(F^{(l)})   $.
\end{itemize}
\end{Definition}

Minimality can be interpreted as a generalization of $L_{2}-$stability. 

\begin{Lemma}\label{lem_min_stab_con}
Let $ F = ( f_{1} , \ldots , f_{N} )^{T}  \in (L_{2}(\mathbb{R}))^{N}  $ be a multiscaling function of multiplicity $N$ according to Definition \ref{def_MRA}, item (vi). If $F$ is $L_{2}-$stable, it is also minimal.
\end{Lemma}

\begin{proof}
Let $ F = ( f_{1} , \ldots , f_{N} )^{T}  \in (L_{2}(\mathbb{R}))^{N}  $ be $L_{2}-$stable. We proceed by contradiction and assume that $F$ is not minimal. Consequently there exists at least one $ l \in \{ 1 , \ldots , N   \}$ such that $f_{l} \in S(F^{(l)})$. With other words there is a representation
\begin{equation}\label{pr_sta_min_eq1}
f_{l}(x) = \sum_{\substack{ n \in \{  1 , \ldots , N   \} \\ n \not = l }} \sum_{k \in \mathbb{Z} } c_{n,k} f_{n} ( x - k ) .
\end{equation}
Using the definition of $L_{2}-$stability, namely Definition \ref{def_L2stab}, we find
\begin{align*}
0 & = \Big \Vert f_{l} -  \sum_{\substack{ n \in \{  1 , \ldots , N   \} \\ n \not = l }} \sum_{k \in \mathbb{Z} } c_{n,k} f_{n} ( \cdot - k ) \Big \vert L_{2}(\mathbb{R}) \Big \Vert^{2} \\
& \geq C_{1} \Big ( 1 + \sum_{\substack{ n \in \{  1 , \ldots , N   \} \\ n \not = l }} \sum_{k \in \mathbb{Z} } | c_{n,k} |^{2}     \Big ) \\
& \geq C_{1} 
\end{align*} 
with $ C_{1} > 0  $. However this is a contradiction. 
\end{proof}

Let us remark that minimality does not imply $L_{2}-$stability, see \cite{Kop1}. Now it is not difficult to see that minimality is valid for our quarklet vector $   \Psi  $ for $ m = 1 $.

\begin{Lemma}\label{lem_gmini1}
Let $ m = 1  $ and $ \tilde{m} \in \mathbb{N}    $ with  $  \tilde{m} + 1 \in 2 \mathbb{N}   $. Let $   p \in \mathbb{N}_{0}  $. Then the set $ \Psi = ( \psi_{0} , \psi_{1} , \ldots , \psi_{p} )^{T} $ is minimal.
\end{Lemma}

\begin{proof}
Let $ m = 1  $ and $ \tilde{m} \in \mathbb{N}    $ with  $  \tilde{m} + 1 \in 2 \mathbb{N}   $. Let $   p \in \mathbb{N}_{0}  $. We already have seen that the integer shifts of the $ p + 1 $ functions
\begin{align*}
\psi_{0} , \psi_{1} , \ldots , \psi_{p} 
\end{align*}
are $L_{2}-$stable, see formula \eqref{eq_proof_polyph_nonmo} in Step 3 of the proof of Lemma \ref{lem_oplus1}. Using Lemma \ref{lem_min_stab_con} we find that those functions are also minimal.  
\end{proof}

Let us remark that it is also possible to prove Lemma \ref{lem_gmini1} by checking Definition \ref{def_gmini}. The proof is then based on the different polynomial degrees of the involved functions. In what follows to apply the theory from \cite{Kaz1} we also need the subsequent definition. 

\begin{Definition}\label{def_per1}
Let $   f , g \in L_{2}(\mathbb{R})  $. Then for all $ x \in \mathbb{R}    $ we put 
\begin{align*}
[ f , g ] (x) := \sum_{k \in \mathbb{Z} } f(x+ 2 \pi k) \overline{g(x+ 2 \pi k)} .
\end{align*}
\end{Definition}

Now we are prepared to provide a describtion for the orthogonal projection onto $ W_{p,0}   $. In the following result we use the convention $ \frac{0}{0} = 0  $.

\begin{Theorem}\label{satz_proj1}
Let $ m = 1  $ and $ \tilde{m} \in \mathbb{N}    $ with $  \tilde{m} + 1 \in 2 \mathbb{N}   $. Let $   p \in \mathbb{N}_{0}  $ and $ \Psi = ( \psi_{0} , \psi_{1} , \ldots , \psi_{p} )^{T} $. By iteration we define $ \psi_{0}^{\star} := \psi_{0}   $ and for $  1 \leq q \leq p   $
\begin{align*}
\mathcal{F} \psi_{q}^{\star} :=  \mathcal{F}  \psi_{q} - \sum_{l=0}^{q-1} \frac{[\mathcal{F}\psi_{q}, \mathcal{F} \psi_{l}^{\star} ]}{[\mathcal{F} \psi_{l}^{\star} , \mathcal{F} \psi_{l}^{\star}]}  \mathcal{F} \psi_{l}^{\star} .
\end{align*}
Then for any $  f \in L_{2}(\mathbb{R})  $ the orthogonal projection $   P_{W_{p,0}}(f)  $ of $ f $ onto $ W_{p,0}   $ is given by
\begin{align*}
\mathcal{F}( P_{W_{p,0}}(f)) = \sum_{q = 0}^{p} \frac{[ \mathcal{F} f , \mathcal{F} \psi_{q}^{\star}]}{[ \mathcal{F} \psi_{q}^{\star}  ,  \mathcal{F} \psi_{q}^{\star}  ]}  \mathcal{F} \psi_{q}^{\star} .
\end{align*}
Moreover for all $ 0 \leq q \leq p  $ the $  [ \mathcal{F} \psi_{q}^{\star}  ,  \mathcal{F} \psi_{q}^{\star}  ]   $ are rational trigonometric functions (see Section 3 in \cite{Kaz1} for a definition).
\end{Theorem}

\begin{proof}
This result is an easy consequence of Lemma 2.1 in \cite{Kaz1}, see also Lemma 2.13 in \cite{TaKaz}.   It is important that the set $ \Psi = ( \psi_{0} , \psi_{1} , \ldots , \psi_{p} )^{T} $ is minimal, see Lemma \ref{lem_gmini1}. The last statement follows from Theorem 3.1 in \cite{Kaz1}. Here we used that the quarklets are compactly supported. 
\end{proof}

The functions $ \psi_{q}^{\star}   $ we obtained in Theorem \ref{satz_proj1} have some pleasant properties. Some of them are collected below. 

\begin{Lemma}\label{lem_prop_ort1}
Let $ m = 1  $ and $ \tilde{m} \in \mathbb{N}    $ with  $ \tilde{m} + 1 \in 2 \mathbb{N}   $. Let $   p \in \mathbb{N}_{0}  $ and $\Psi^{\star} = ( \psi_{0}^{\star} , \psi_{1}^{\star} , \ldots , \psi_{p}^{\star} )^{T} \in ( L_{2}(\mathbb{R}) )^{p+1} $ be the functions defined in Theorem \ref{satz_proj1}. Then the following assertions are true. 
\begin{itemize}
\item[(i)] For all $  0 \leq q \leq p  $ we have $ \psi_{q}^{\star} \in  W_{p,0}    $.
\item[(ii)] For all $  0 \leq q , r \leq p    $ with $  q \not = r  $ we observe $   \langle \psi_{q}^{\star}  ,  \psi_{r}^{\star}    \rangle_{L_{2}(\mathbb{R})} = 0     $.
\end{itemize}
\end{Lemma}

\begin{proof}
\textit{Step 1.} At first we prove (i).  This follows directly from Proposition 2.2 in \cite{Kaz1}.  We use that the set $ \Psi = ( \psi_{0} , \psi_{1} , \ldots , \psi_{p} )^{T} $ is minimal, see Lemma \ref{lem_gmini1}. Recall that we have $ W_{p,0}  = S (\Psi)    $.

\textit{Step 2.} Now we prove (ii).  For that purpose at first we use the Plancherel Theorem. Then we find
\begin{align*}
  \langle \psi_{q}^{\star}  ,  \psi_{r}^{\star}    \rangle_{L_{2}(\mathbb{R})}  & = \int_{\mathbb{R}}  \psi_{q}^{\star}(x) \overline{ \psi_{r}^{\star}(x)} dx \\
& = \int_{\mathbb{R}} (\mathcal{F} \psi_{q}^{\star})(\xi) \overline{ (\mathcal{F} \psi_{r}^{\star})(\xi)} d \xi \\ 
& = \sum_{k \in \mathbb{Z}} \int_{2 \pi k}^{2 \pi k+ 2 \pi} (\mathcal{F} \psi_{q}^{\star})(\xi) \overline{ (\mathcal{F} \psi_{r}^{\star})(\xi)} d \xi \\ 
& =  \int_{0}^{2 \pi} \sum_{k \in \mathbb{Z}}  (\mathcal{F} \psi_{q}^{\star})(\xi + 2 \pi k) \overline{ (\mathcal{F} \psi_{r}^{\star})(\xi + 2 \pi k)} d \xi .
\end{align*}
Next we apply Definition \ref{def_per1}. We obtain
\begin{align*}
  \langle \psi_{q}^{\star}  ,  \psi_{r}^{\star}    \rangle_{L_{2}(\mathbb{R})}  & =  \int_{0}^{2 \pi} [ \mathcal{F} \psi_{q}^{\star} ,   \mathcal{F} \psi_{r}^{\star}  ](\xi) d \xi .
\end{align*}
Formula (2.4) in Proposition 2.2 in \cite{Kaz1} yields $  [ \mathcal{F} \psi_{q}^{\star} ,   \mathcal{F} \psi_{r}^{\star}  ](\xi)  = 0     $ almost everywhere on $   [0,2 \pi]  $. Consequently we get $ \langle \psi_{q}^{\star}  ,  \psi_{r}^{\star}    \rangle_{L_{2}(\mathbb{R})}   =  0 $ . The proof is complete.
\end{proof}

Let us give an example to illustrate how the functions $ \psi_{q}^{\star}   $ look like. 

\begin{example}\label{ex_mod_quarklets}
Let $ m = \tilde{m} = 1 $ and $ p = 3 $. That means we work with $\Psi = ( \psi_{0} , \psi_{1} , \psi_{2} , \psi_{3} )^{T} $. Then for $ \Psi^{\star}  $ we find
\begin{align*}
        \psi_0^\star &= \psi_0 =  \begin{cases*}
        1, & $x \in [0,1/2),$ \\
        -1, &  $x \in [1/2,1), $
        \end{cases*}\\
        \psi_1^\star &= \psi_1 - \frac12 \psi_0^\star =  \begin{cases*}
    2x-\frac12, & $x \in [0,1/2),$ \\
        -2x+ \frac32, &  $x \in [1/2,1), $
        \end{cases*}\\
        \psi_2^\star &=\psi_2 - \psi_1^\star - \frac13 \psi_0^\star =
\psi_2-\psi_1+\frac16\psi_0 =  \begin{cases*}
        4x^2-2x+\frac16, & $x \in [0,1/2),$ \\
        -4x^2+6x-\frac{13}{6}, &  $x \in [1/2,1). $
        \end{cases*}
       \end{align*}
For $q = 3$ we have
\begin{align*}
 \psi_3^\star &= \psi_3 -\frac32 \psi_2 + \frac35 \psi_1-\frac{1}{20}\psi_0 =
\begin{cases*}
        8x^3-6x^2+\frac65 x -\frac{1}{20}, & $x \in [0,1/2),$ \\
        -8x^3+18x^2-\frac{66}{5}x+\frac{63}{20}, &  $x \in [1/2,1). $
        \end{cases*}    
\end{align*}
Notice that the branches of the resulting functions remind on the Legendre polynomials transplanted onto the subintervals.

\begin{figure}[h]

\begin{minipage}[b]{.4\linewidth}

\includegraphics[width=\linewidth]{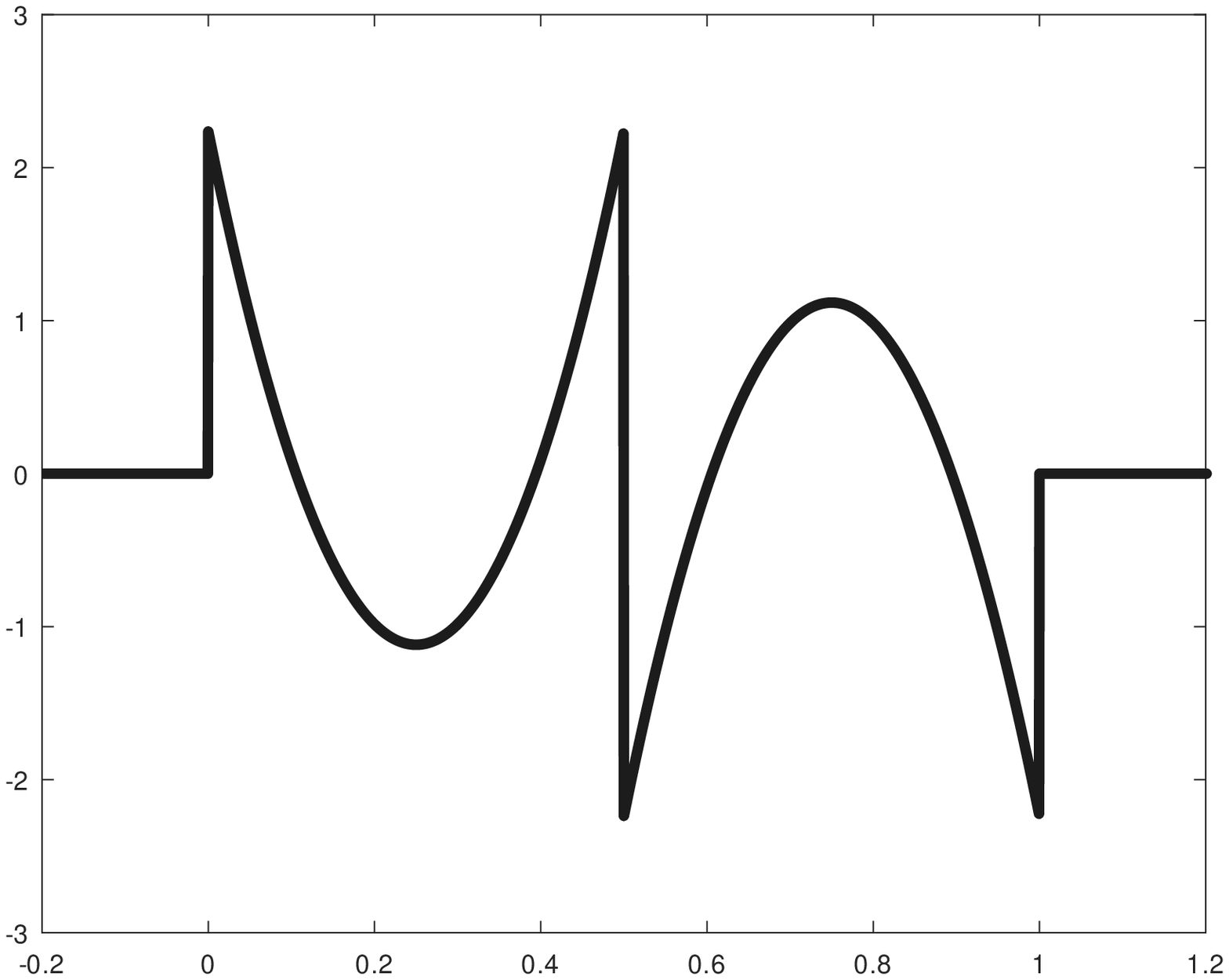}

\caption{$  \psi_2^\star  $ for $ m = \tilde{m} = 1 $.}

\end{minipage}
\hspace{.1\linewidth}
\begin{minipage}[b]{.4\linewidth}

\includegraphics[width=6cm]{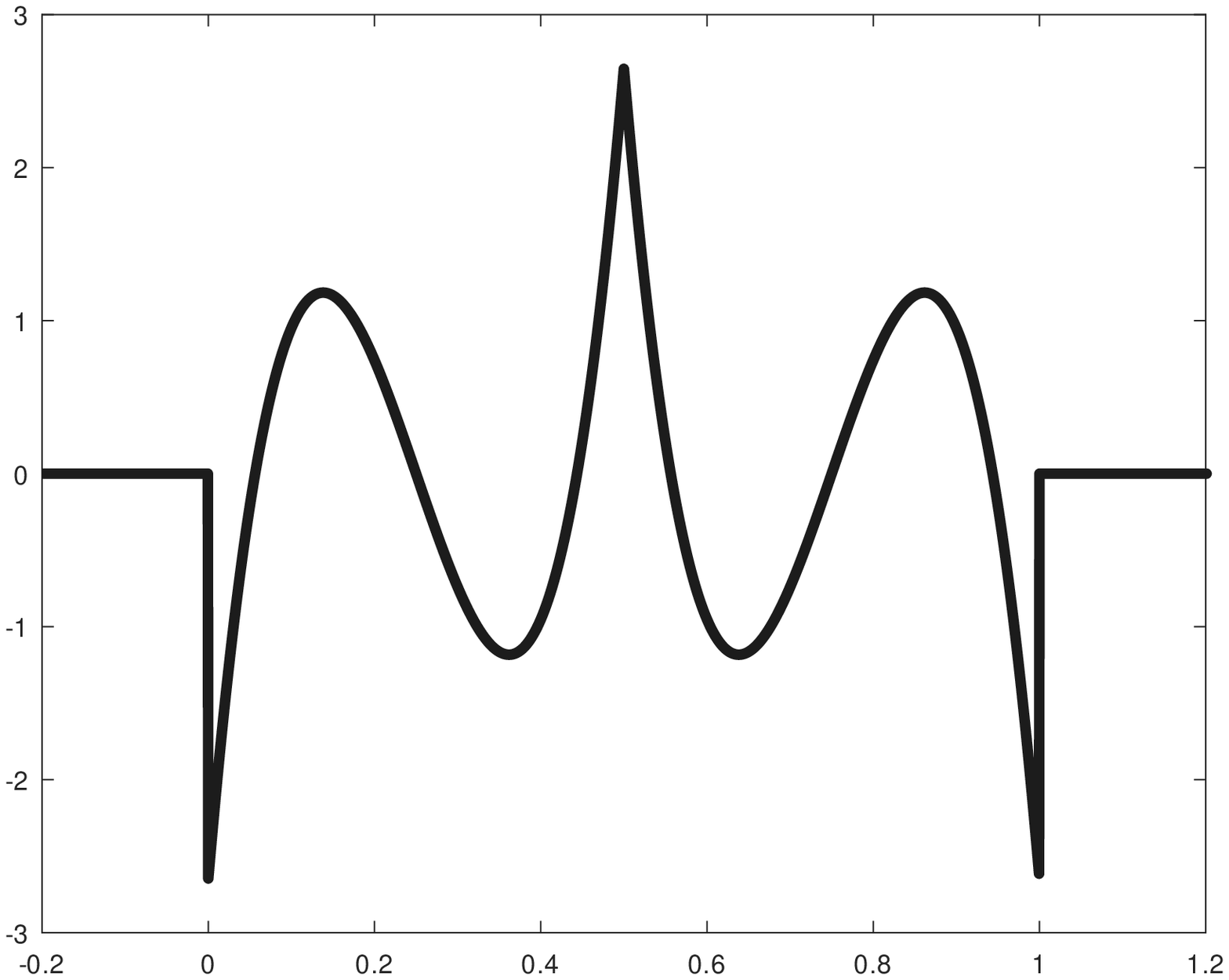}

\caption{$  \psi_3^\star  $ for $ m = \tilde{m} = 1 $.}

\end{minipage}

\end{figure}

\end{example}

\begin{Remark}
Notice that for $ m = 1  $ the modified quarklets collected in $ \Psi^{\star}  $ are linear combinations of the original quarklets $ \Psi  $, see also Example \ref{ex_mod_quarklets}. Therefore the functions $ \psi_{q}^{\star} $ possess at least the same number of vanishing moments as the original quarklets $ \psi_{q} $. For results concerning vanishing moments of quarklets  we refer to \cite{bib:DKR17}.
\end{Remark}

We can use the modified quarklets from $ \Psi^{\star}    $ to rewrite the spaces $ W_{p,j}     $ in a more convenient way. For that purpose let us give the following definition.

\begin{Definition}\label{def_Wpj_star}
Let $ m = 1  $ and $ \tilde{m} \in \mathbb{N}    $ with  $  \tilde{m} + 1 \in 2 \mathbb{N}   $. Let $  j \in \mathbb{Z}  $ and $   p \in \mathbb{N}_{0}  $. Let $ \Psi^{\star} = ( \psi_{0}^{\star} , \psi_{1}^{\star} , \ldots , \psi_{p}^{\star} )^{T} $ be the functions defined in Theorem \ref{satz_proj1}. Then we define
\begin{align*}
W_{p,j}^{\star} :=  \overline{ \spa \{ 2^{\frac{j}{2}} \psi_{q}^{\star} (2^{j} \cdot - k ) : 0 \leq q \leq p , k \in \mathbb{Z}   \} } \subset L_{2}(\mathbb{R}).
\end{align*}
Moreover for each $  0 \leq q \leq p  $ we put
\begin{align*}
W_{q,j}^{\star \star} :=  \overline{ \spa \{ 2^{\frac{j}{2}} \psi_{q}^{\star} (2^{j} \cdot - k ) :  k \in \mathbb{Z}   \} } \subset L_{2}(\mathbb{R}).
\end{align*}
\end{Definition}

It is not difficult to see that the spaces $ W_{p,j}  $ and $   W_{p,j}^{\star} $ are equal. 

\begin{Lemma}\label{lem_orth_zerl1}
Let $ m = 1  $ and $ \tilde{m} \in \mathbb{N}    $ with $  \tilde{m} + 1 \in 2 \mathbb{N}   $. Let $   j \in \mathbb{Z}    $ and $   p \in \mathbb{N}_{0}  $. Then we have
\begin{align*}
W_{p,j}^{\star} = \hat{\bigoplus}_{q = 0}^{p}  W_{q,j}^{\star \star} = W_{p,j} .
\end{align*}
Here $  \hat{\oplus}  $ stands for the orthogonal direct sum.
\end{Lemma}

\begin{proof}
To prove this result we use Proposition 2.3 from \cite{Kaz1}. Recall that we have $  S(\Psi^{\star})  =  W_{p,0}^{\star}  $ and $    S( \psi_{q}^{\star} ) =  W_{q,0}^{\star \star}  $. The set $  \Psi $ is minimal, see Lemma \ref{lem_gmini1}. Consequently the assertion of Lemma \ref{lem_orth_zerl1} follows when we apply Proposition 2.3 from \cite{Kaz1}. 
\end{proof}

Using the spaces from Definition \ref{def_Wpj_star} we can reformulate Theorem \ref{satz_L2zerl1}.

\begin{Theorem}\label{satz_L2zerl2}
Let $ m = 1  $ and $ \tilde{m} \in \mathbb{N}    $ with  $  \tilde{m} + 1 \in 2 \mathbb{N}   $. Let $   p \in \mathbb{N}_{0}  $. Then we have
\begin{align*}
L_{2}(\mathbb{R} ) = \overline{ \bigoplus_{j \in \mathbb{Z}} W_{p,j}^{\star} } = \overline{  \bigoplus_{j \in \mathbb{Z}} \hat{\bigoplus}_{q = 0}^{p}  W_{q,j}^{\star \star} }
\end{align*}
and
\begin{align*}
L_{2}(\mathbb{R} ) = \overline{  V_{p,0} \oplus \bigoplus_{j = 0}^{\infty}  W_{p,j}^{\star} }  = \overline{ V_{p,0} \oplus \bigoplus_{j = 0}^{\infty}  \hat{\bigoplus}_{q = 0}^{p}  W_{q,j}^{\star \star} } .
\end{align*}
Here $  \oplus  $ denotes the (non-orthogonal) algebraic direct sum. Moreover $  \hat{\oplus}  $ stands for the orthogonal direct sum.
\end{Theorem}

\begin{proof}
To prove this result on the one hand we use Theorem \ref{satz_L2zerl1}. There we obtained
\begin{align*}
L_{2}(\mathbb{R} ) = \overline{ \bigoplus_{j \in \mathbb{Z}} W_{p,j} } = \overline{ V_{p,0} \oplus \bigoplus_{j = 0}^{\infty} W_{p,j} } .
\end{align*}
On the other hand we apply Lemma \ref{lem_orth_zerl1}. Now a combination of both observations completes the proof. 
\end{proof}

\begin{Remark}
Let $m=1$ and $ p \in \mathbb{N} $. Then there exists an alternative way to construct a vector $ \breve{\Phi} \in  (L_{2}(\mathbb{R}))^{p+1}   $ which is biorthogonal to the quark vector $ \Phi =  ( \varphi_{0} , \varphi_{1} , \ldots , \varphi_{p}    )^{T}   $ in the sense of Definition \ref{def_vec_biorto}. For that purpose we define the Gramian matrix $G$ which is given by
\begin{align*}
G := \big ( [ \mathcal{F} \varphi_{i-1}  , \mathcal{F} \varphi_{j-1} ]   \big )_{i=1, \ldots , p+1 ; j = 1, \ldots , p+1} .
\end{align*} 
Since for $m=1$ the quark vector $ \Phi $ is $L_{2}-$stable, see Lemma \ref{lem_L2stab_multiquark}, the existence of $ G^{-1}  $ follows. For that we refer to \cite{JetPlo}, see Section 2.2.9. Now we introduce a function vector $ \breve{\Phi} $ given by
\begin{align*}
\breve{\Phi} := \mathcal{F}^{-1}[G^{-1} \mathcal{F} \Phi] .
\end{align*} 
It has been shown in \cite{Nil}, see Proposition 2.1, that $ \Phi  $ and $ \breve{\Phi}  $ are biorthogonal in the sense of Definition \ref{def_vec_biorto}. Here we also can refer to Theorem 2.2.14 in \cite{JetPlo}. An advantage of this approach is that we have $ \breve{\Phi} \in (  L_{2}(\mathbb{R})  )^{p+1}  $ and that $ \Phi  $ and $ \breve{\Phi}   $ are dual in the usual $L_{2}-$sense. Moreover using the usual procedure now it is possible to construct a biorthogonal quarklet system. However the resulting quarklet functions do not coincide with $  \Psi  $ and are not compactly supported, but they are decaying exponentially at least.
\end{Remark}

\section{Appendix}\label{sec_appendix}

Here we present the proof of Lemma \ref{lem_stab_sinleres}. It is rather technical since we have to compute the zeros of several functions. 

\textit{Proof of Lemma \ref{lem_stab_sinleres}. }

\textit{Step 1.}
To prove this result we use Lemma \ref{lem_stab_singl}. Notice that the Fourier transforms of the quarks $ \mathcal{F} \varphi_{l}    $ already have been calculated in Lemma \ref{lem_qua_fourtrans}. Now we explicitly state the functions $   \mathcal{F} \varphi_{l}   $ for several $ m \in \mathbb{N} $ and $  l \in \mathbb{N}_{0}    $. Moreover we calculate some of their zeros. For each $ m \in \mathbb{N}   $ and $ l=0 $ the function $  \mathcal{F} \varphi_{0}   $ has zeros $ \xi_{k} = 2 k \pi    $ with $ k \in \mathbb{Z}   \setminus \{ 0 \} $. For $m=1 $ and $ l=0$ we find
\[
( \mathcal{F} \varphi_{0} ) (\xi) = \left\{ \begin{array}{lll}    e^{- i \xi \frac{1}{2}}\frac{\sin(\frac{\xi}{2})}{\frac{\xi}{2}}   & \qquad & \mbox{for} \qquad \xi \not = 0   ;
\\  
1 & \qquad & \mbox{for}\qquad  \xi = 0 .
\\
\end{array}
\right.
\]
This function has the zeros $ \xi_{k} = 2 k \pi    $ with $ k \in \mathbb{Z}   \setminus \{ 0 \} $. For $m=1 $ and $ l=1$ we find
\[
( \mathcal{F} \varphi_{1} ) (\xi) = \left\{ \begin{array}{lll}    e^{- i \xi \frac{1}{2}} \Big ( \frac{ \xi \cos(\frac{\xi}{2}) + (-2- i \xi)  \sin(\frac{\xi}{2})}{\xi^{2}} \Big )  & \qquad & \mbox{for} \qquad \xi \not = 0   ;
\\  
- \frac{1}{2} i & \qquad & \mbox{for}\qquad  \xi = 0 .
\\
\end{array}
\right.
\]
One can observe that this function has no zero in $ [-10 , 10 ]     $.  For $m=1 $ and $ l=2$ we have
\[
( \mathcal{F} \varphi_{2} ) (\xi) = \left\{ \begin{array}{lll}    e^{- i \xi \frac{1}{2}} \Big ( \frac{ (-2 - i \xi ) \xi \cos(\frac{\xi}{2}) + ( - \xi^{2} + 2 i \xi + 4)  \sin(\frac{\xi}{2})}{\xi^{3}} \Big )  & \qquad & \mbox{for} \qquad \xi \not = 0   ;
\\  
- \frac{1}{3}  & \qquad & \mbox{for}\qquad  \xi = 0 .
\\
\end{array}
\right.
\]
This function has no zero in $  [-7 , 7 ]   $.  For $m=1 $ and $ l=3$ we have
\[
( \mathcal{F} \varphi_{3} ) (\xi) = \left\{ \begin{array}{lll}    e^{- i \xi \frac{1}{2}} \Big ( \frac{ ( - \xi^{2} + 3 i \xi + 6 ) \xi \cos(\frac{\xi}{2}) + ( i \xi^{3} + 3 \xi^{2} - 6 i \xi -12 )  \sin(\frac{\xi}{2})}{\xi^{4}} \Big )  & \qquad & \mbox{for} \qquad \xi \not = 0   ;
\\  
 \frac{1}{4}  i & \qquad & \mbox{for}\qquad  \xi = 0 .
\\
\end{array}
\right.
\]
There are no zeros in $  [ -7 , 7 ]    $. For $ m $ even and $ l = 0   $ we find
\[
( \mathcal{F} \varphi_{0} ) (\xi) = \left\{ \begin{array}{lll}  \Big ( \frac{\sin(\frac{\xi}{2})}{\frac{\xi}{2}} \Big ) ^{m}  & \qquad & \mbox{for} \qquad \xi \not = 0   ;
\\  
1 & \qquad & \mbox{for}\qquad  \xi = 0 .
\\
\end{array}
\right.
\]
This function has the zeros $ \xi_{k} = 2 k \pi    $ with $  k \in \mathbb{Z} \setminus \{ 0 \}   $. For $ m $ even and $ l = 1   $ we get 
\[
( \mathcal{F} \varphi_{1} ) (\xi) = \left\{ \begin{array}{lll} m \Big ( \frac{\sin(\frac{\xi}{2})}{\frac{\xi}{2}} \Big ) ^{m-1} \Big ( \frac{\cos(\frac{\xi}{2})}{\xi} - \frac{2 \sin (\frac{\xi}{2})}{\xi^{2}}  \Big )  & \qquad & \mbox{for} \qquad \xi \not = 0   ;
\\  
0 & \qquad & \mbox{for}\qquad  \xi = 0 .
\\
\end{array}
\right.
\]
This function has (beside others) periodic zeros of the form $ \xi_{k} = 2 k \pi     $ with $ k \in \mathbb{Z}    $. For $ m = 2  $ and $ l = 2   $ we observe 
\[
( \mathcal{F} \varphi_{2} ) (\xi) = \left\{ \begin{array}{lll} \frac{2( - (\xi^{2} - 12) \sin^{2}(\frac{\xi}{2}) + \xi^{2} \cos^{2}(\frac{\xi}{2}) - 8 \xi \sin (\frac{\xi}{2 })  \cos (\frac{\xi}{2}))  }{\xi^{4}}  & \qquad & \mbox{for} \qquad \xi \not = 0   ;
\\  
- \frac{1}{6} & \qquad & \mbox{for}\qquad  \xi = 0 .
\\
\end{array}
\right.
\]
This function has zeros of the form $ \xi_{1,2} \approx \pm 2,606       $, $    \xi_{3,4} \approx \pm 7,414$, $ \xi_{5,6} \approx \pm 10,562      $ and so forth. For $ m = 2    $ and $ l = 3 $ the function $ \mathcal{F} \varphi_{3}  $ has in $ [-2 \pi , 2 \pi ]     $ the zeros $ \xi_{1} = 0   $ and $ \xi_{2,3} \approx \pm 4,639  $.  For $m=3 $ and $ l=0$ we find
\[
( \mathcal{F} \varphi_{0} ) (\xi) = \left\{ \begin{array}{lll}    e^{- i \xi \frac{1}{2}} \Big ( \frac{\sin(\frac{\xi}{2})}{\frac{\xi}{2}} \Big )^{3}  & \qquad & \mbox{for} \qquad \xi \not = 0   ;
\\  
1 & \qquad & \mbox{for}\qquad  \xi = 0 .
\\
\end{array}
\right.
\]
This function has the zeros $ \xi_{k} = 2 k \pi    $ with $  k \in \mathbb{Z} \setminus \{ 0 \}   $. For $m=3 $ and $ l=1$ we have
\[
( \mathcal{F} \varphi_{1} ) (\xi) = \left\{ \begin{array}{lll}    e^{- i \xi \frac{1}{2}} \sin^{2}(\frac{\xi}{2}) \frac{( 12 \xi \cos (\frac{\xi}{2}) + ( -24 - 4 i \xi ) \sin (\frac{\xi}{2}))}{\xi^{4}}  & \qquad & \mbox{for} \qquad \xi \not = 0   ;
\\  
- \frac{1}{2} i & \qquad & \mbox{for}\qquad  \xi = 0 .
\\
\end{array}
\right.
\]
In $ [ - 2 \pi , 2 \pi ]     $ this function only has the zeros $ \xi_{1,2} = \pm 2 \pi     $.  For $m=3 $ and $ l=2$ we have
\[
( \mathcal{F} \varphi_{2} ) (\xi) = \left\{ \begin{array}{lll}    e^{- i \xi \frac{1}{2}} \sin(\frac{\xi}{2}) \frac{( - 8 \xi^{2} + 24 i \xi + 96 ) \sin^{2}(\frac{\xi}{2}) + 12 \xi^{2} \cos^{2}(\frac{\xi}{2}) + ( - 72 - 12 i \xi ) \xi \sin(\frac{\xi}{2}) \cos(\frac{\xi}{2}) }{\xi^{5}}  & \quad  \xi \not = 0   ;
\\  
- \frac{1}{2}  & \quad  \xi = 0 .
\\
\end{array}
\right.
\]
In $ [ - 2 \pi , 2 \pi ]     $ this function only has the zeros $ \xi_{1,2} = \pm 2 \pi     $. For $ m = 3    $ and $  l = 3   $ the function $ \mathcal{F} \varphi_{3}   $ has no zeros in $ [ - 6 , 6 ]     $.  For $ m = 4  $ and $ l = 2   $ we observe 
\[
( \mathcal{F} \varphi_{2} ) (\xi) = \left\{ \begin{array}{lll} - \frac{ 16 \sin^{2}(\frac{\xi}{2})(( \xi^{2} - 20 ) \sin^{2}(\frac{\xi}{2}) - 3 \xi^{2} \cos^{2}(\frac{\xi}{2}) + 16 \xi \sin (\frac{\xi}{2 })  \cos (\frac{\xi}{2}))  }{\xi^{6}}  & \qquad & \mbox{for} \qquad \xi \not = 0   ;
\\  
- \frac{1}{3} & \qquad & \mbox{for}\qquad  \xi = 0 .
\\
\end{array}
\right.
\] 
This function has zeros of the form $ \xi_{1,2} \approx \pm 1,780       $, $    \xi_{3,4} \approx \pm 7,938$, $ \xi_{5,6} \approx \pm 10,036      $ and so forth. For $ m = 4   $ and $ l = 3    $ the function $ \mathcal{F} \varphi_{3}    $  has (beside others) periodic zeros of the form $ \xi_{k} = 2 k \pi     $ with $ k \in \mathbb{Z}    $. For $m$ odd with $ m > 1  $ and $ l=1 $ we find
\[
( \mathcal{F} \varphi_{1} ) (\xi) = \left\{ \begin{array}{lll}    e^{- i \xi \frac{1}{2}} \Big ( \frac{ \sin(\frac{\xi}{2})} {\frac{\xi}{2}} \Big )^{m-1} \Big (  - \frac{i}{2} \frac{ \sin(\frac{\xi}{2})} {\frac{\xi}{2}} + m \frac{\frac{1}{4} \xi \cos(\frac{\xi}{2}) - \frac{1}{2} \sin(\frac{\xi}{2}) }{(\frac{\xi}{2})^2}  \Big )  & \qquad & \mbox{for} \qquad \xi \not = 0   ;
\\  
- \frac{1}{2} i & \qquad & \mbox{for}\qquad  \xi = 0 .
\\
\end{array}
\right.
\]
In $ [ - 2 \pi , 2 \pi ]    $ this function only has the zeros $  \xi_{1,2} = \pm 2 \pi    $.

\textit{Step 2.} Now we apply Lemma \ref{lem_stab_singl}. Since the quarks are compactly supported functions from $ L_{2}(\mathbb{R}) $ this is possible. We use different combinations of $ m \in \mathbb{N} $ and $ l \in \mathbb{N}_{0} $ and check whether we have 
\begin{equation}\label{eq_stab_prov1}
\sum_{\alpha \in \mathbb{Z}} | (\mathcal{F} \varphi_{l} ) (\xi + 2 \pi \alpha )    |^{2} > 0 
\end{equation} 
for all $ \xi \in \mathbb{R}  $ or not. For that purpose we use the results from Step 1. If we find a $ 2 \pi$ periodic zero then \eqref{eq_stab_prov1} cannot be fulfilled for all $ \xi \in \mathbb{R}  $. If there is no $ 2 \pi$ periodic zero then \eqref{eq_stab_prov1} holds. So we can use Lemma \ref{lem_stab_singl} to complete the proof. \hfill $ \qed $


\vspace{1 cm}

\textbf{Funding.} This paper is a result of the DFG project 'Adaptive high-order quarklet frame methods for elliptic operator equations` with grant numbers DA$360/24-1$ (Marc Hovemann) and RA$2090/3-1$ (Thorsten Raasch).

\vspace{0,5 cm}

\textbf{Acknowledgment.} The authors would like to thank Stephan Dahlke for several tips and hints.

\end{document}